\documentclass[utf8,review]{elsarticle}

\usepackage{lineno,hyperref}
\modulolinenumbers[5]

\journal{Journal of \LaTeX\ Templates}

\usepackage{fullpage}
\usepackage{amsmath}
\usepackage{mathrsfs}
\usepackage{amsfonts}
\usepackage{amssymb}
\usepackage{bm}
\usepackage{enumerate}
\usepackage{tikz}
\usepackage{graphicx}
\usepackage{stfloats}
\usepackage{multicol}

\def\i{{\bf i}}

\def \A {\mathbf{A}}

\newtheorem{lemma}{Lemma}
\newtheorem{theorem}{Theorem}
\newtheorem{remark}{Remark}

\newtheorem{definition}{Definition}
\newtheorem{example}{Example}
\newtheorem{corollary}{Corollary}
\newenvironment{nproof}{\noindent{\em \textbf{Proof.}}}{\quad \hfill$\Box$\vspace{2ex}}

\begin{document}

\begin{frontmatter}


\title{Convolution theorems associated with quaternion linear canonical transform and applications}
\author{Xiaoxiao Hu}
\ead{huxiaoxiao@wmu.edu.cn}
\address{The First Affiliated Hospital of Wenzhou Medical University, Wenzhou Medical University, Wenzhou, Zhejiang, China}

\author{Dong Cheng}
\ead{chengdong720@163.com}
\address{Research Center for Mathematics and Mathematics Education, Beijing Normal University at Zhuhai, Zhuhai, China}
\author{Kit Ian Kou \corref{mycorrespondingauthor}} 

\cortext[mycorrespondingauthor]{Corresponding author}
\ead{kikou@umac.mo}

\address{Department of Mathematics, Faculty of Science and Technology, University of Macau, Macao, China}

%
%
%

\begin{abstract}
Novel types of convolution operators for quaternion linear canonical transform (QLCT) are proposed. Type one and two are defined in the spatial and QLCT spectral domains, respectively. They are distinct in the quaternion space and are consistent once in complex or real space. Various types of convolution formulas are discussed. Consequently, the QLCT of the convolution of two quaternionic functions can be implemented by the product of their QLCTs, or the summation of the products of their QLCTs. As applications, correlation operators and theorems of the QLCT are derived. The proposed convolution formulas are used to solve Fredholm integral equations with special kernels. Some systems of second-order partial differential equations, which can be transformed into the second-order quaternion partial differential equations, can be solved by the convolution formulas as well. As a final point, we demonstrate that the convolution theorem facilitates the design of multiplicative filters.
\end{abstract}


\end{frontmatter}

\section{Introduction}\label{sec1}
The convolution theorem plays a crucial role in signal processing \cite{wei2020convolution} and solving differential and integral equations \cite{PWJ2000}. It states that under suitable conditions the Fourier transform of a convolution equals the product of their Fourier transforms.  It would be convenient for the analysis of multiplicative filters in the spectral domain.  The Fourier transform is a special case of linear canonical transform, both of them are useful tools in signal processing. The convolution theorem of the linear canonical transform has been investigated  \cite{wei2020convolution,deng2010comments,shi2014generalized,deng2006convolution,deng2010comments,wei2012convolution,wei2009convolution}. These convolution structures have elegance and simplicity, which states that a modified ordinary convolution in the time domain is equivalent to a simple multiplication operation for linear canonical transform and Fourier transform. Nevertheless, these convolution structures are only suitable for real and complex signals, not for higher dimesional signals like quaternions. Quaternions have shown advantages over real and complex within multidimensional signal processing due to the good consistency of the correlation among its four components \cite{kou2013uncertainty,kou2014asymptotic}.
A  color image can be expressed by a  pure quaternion number. As opposed to processing each of the three components independently, all three color components of the image are treated together. Multidimensional signal processing by quaternions has become a hot research topic \cite{wang2021robust}. 
The quaternion Fourier transform (QFT) is a powerful tool in colour image processing \cite{ell2013quaternion,le2014instantaneous}.
To meet the new needs of quaternion signal processing, the quaternion linear canonical transform (QLCT) is developed to circumvent this problem. A linear canonical transform is a generalization of the Fourier transform, so researchers introduced the QLCT \cite{kou2013uncertainty}, which includes many well-known signal processing operations,  such as the QFT and quaternion fractional Fourier transform \cite{guanlei2008fractional, serbes2010optimum}  as its special cases.  Furthermore, the quaternionic offset linear canonical transform is a generalization of the QLCT and has been studied in  \cite{bhat2022wvd,bhat2022algebra}. The QLCT was shown to be a flexible tool in signal processing \cite{guo2011reduced,hu2021sampling, chen2015pitt,kou2014asymptotic,yang2014uncertainty,kou2016envelope}.  Some important theorems of QLCT have been studied, such as the inversion theorems \cite{hu2016quaternion},  the Phancherel theorems \cite{kou2019Plancherel}, uncertainly principles \cite{kou2013uncertainty,yang2014uncertainty}, and sampling theorems\cite{hu2021sampling,hu2022sampling}.
 
 \subsection{Related works}
The convolution theorem associated with the quaternion Fourier transform (QFT) has been studied \cite{bahri2013convolution, bujack2014convolution, hitzer2016general, de2015connecting}. The convolution theorems for the two-sided QFT are studied in \cite{bahri2013convolution}, which stated the  QFT of the convolution operator is equal to the summation of the products of their QFTs. Bahri et al. applied their QFT convolution to study hypoelliptic and to solve the heat equation in \cite{bahri2013convolution}. Later Bujack et al. proposed the Mustard convolution of the two-sied QFT in \cite{bujack2014convolution}. The Mustard convolution can be expressed in the spectral domain as the pointwise product of the QFTs of the factor functions. It behaves nicely in the spectral domain but complicates in the corresponding spatial domain. Hitzer et al. studied the  Mustard convolution of the two-sied QFT in \cite{ hitzer2016general}.  The relationship between the classical and the Mustard convolutions is studied in detail \cite{ hitzer2016general}.   De Bie et al.  proposed a new Mustard convolution for the left-sided QFT, the correspondence between classical and the new Mustard convolution is analyzed in \cite{de2015connecting}.

Recently, researchers introduced the QLCT convolution theorems \cite{sam2021new, bahri2019Two, li2021new}. 
Saima et al. in \cite{sam2021new} studied the theory of one-dimensional QLCT, such as the inversion formula, linearity, convolution theorem, Parseval's identity, and the product theorem.  There are mainly two types of QLCT in two-dimensional QLCT,   one-sided and two-sided QLCT, due to the incommutable nature of quaternion multiplication. 
 The one-sided QLCT was proposed by Mawardi et al.by in \cite{bahri2019Two}. They studied its convolution and correlation theorems and established its uncertainty principle.
Li et al .  introduced a canonical convolution operator for the two-sided QLCT in \cite{li2021new}. They studied its convolution, correlation, and product theorems.
In contrast to the Fourier transform, they do not have the elegant classical result or efficient implementation. The convolution theorem derived in Li et al. \cite{li2021new} shows that the two-sided QLCT of the convolution of the two quaternionic functions is equal to the sum of the products of their QLCTs and two-phase factors. Hence it is not the same as the Fourier transform theorem. It is not easy to implement comparable to that of the Fourier transform.

In the present paper, two novel types of convolution operators for the two-sided QLCT are proposed. They are defined on spatial and spectral domains respectively. They are distinct in the quaternion space and are consistent once in complex or real space.
These are different from those introduced in \cite{li2021new, bahri2019Two}. The proposed convolution theorems state that the QLCT of convolution in the spatial domain is equivalent to a multiplication operation of their QLCTs or summation of multiplication operations. Using the new convolutions, we define the correlation operators,  solve the Fredholm integral equation, some special systems of second-order partial differential equations, and multiplication filters in QLCT.
 
 \subsection{Paper Contributions}
The contributions of this paper are summarized as
follows.
\begin{enumerate}
	\item Two types of QLCT convolution operators are proposed, denoted by  $\circledast^{\bm{\mu},\bm{\nu}}_{A_{1},A_{2}}$ and  $ \bigstar$. Their definitions are based on spatial and spectral dimensions, respectively. The corresponding convolution theorems,  the Parseval formula, the energy theorem, and the product theorem are established.
	\begin{itemize}
		\item The convolution theorem for quaternion slice functions is given by  Theorem  \ref{le43}.
		\item The convolution theorems for  quaternion   functions are provided in  Theorems \ref{h1}, \ref{th46} and \ref{h0312}.
	\end{itemize}
	
	\item Four types of applications of convolution theorems are given.
	\begin{itemize}
	 \item The first one is to define two correlation operators 
	by the corresponding convolution operators, which can be reduced to the correlation operations for 2-D Fourier transform and QFT in  \cite{ell2014quaternion}.  Moreover, the correlation theorems are obtained.
	  \begin{itemize}
		\item The correlation theorem for slice quaternion functions is given by Theorem    \ref{th62}.
		\item The correlation theorems for quaternion   functions are provided in  Theorems  \ref{th63} and  \ref{th262}.
	   \end{itemize}
	\item The second one is to solve the  Fredholm integral equation of the first kind involving special kernels.
   \item The third one is to solve some special system of second-order partial differential equations. 
  \item The last one is to design multiplication filters in the QLCT domain.
	\end{itemize}
\end{enumerate}
\subsection{Paper Outlines}
The paper is organized as follows.  Section \ref{covsec2} gives a brief review of
quaternion,  QFT,  and QLCT as well as their basic properties.
Section \ref{covsec3}   proposes the first type of convolution structure for the QLCT and applications.
Section \ref{covsec4}  provides the second type of convolution structure for the QLCT and applications.
Some conclusions are drawn, and future works are proposed in Section \ref{covsec5}.

\subsection{Notations}
In this work, scalars are represented using an italic letter (e.g., $x$). The three imagery units of quaternion space are denoted by
boldface lowercase letters (e.g., ${\bm \mu}$, ${\bm \nu}$, and ${\bm \eta}$). Table \ref{tab:notation} summarizes the key notations and acronyms used in this paper.

\begin{table}
\caption{List of acronyms and symbols used in this article.}
\begin{center}
\begin{tabular}{l|l}
\hline\cline{1-2}
\textbf{Acronym}       &\textbf{Description}
\\ \hline
$\mathbb{H}$  & Quaternion space.   \\
 $ \mathbb{U}$ &  Unit pure imaginary  quaternion, $ \mathbb{U}:=\{ q| q=q_r+\underline{q}, q_r=0, \underline{q}^{2}=-1\}.$  \\
$\mathbb{H}(\bm{\mu})$  &  Slice quaternion, $\mathbb{H}(\bm{\mu}):=\{ q | q=q=q_{r}+\bm{\mu}q_{\mu}, q_{r}, q_{\mu} \in \mathbb{R}, \bm{\mu} \in \mathbb{U}\}.$   \\
 $L^{p}(\mathbb{R}^2, \mathbb{H}) $ & $\{ f| f:\mathbb{R}^2 \to \mathbb{H}, \Vert f \Vert_{p}:=\left (\int_{\mathbb{R}^2} |f(x,y)|^p dxdy \right )^{\frac{1}{p}} <\infty \}$.\\
   $\mathcal{L}^{\bm{\mu},\bm{\nu}}[f]$ &  $\mathcal{L}^{\bm{\mu},\bm{\nu}}[f] =\mathcal{L}^{\bm{\mu},\bm{\nu}}_{A_{1},A_{2}}[f]$.\\
 $A_{i}$  &   $A_{i}=\left(
             \begin{array}{cc}
               a_{i} &b_{i} \\
               c_{i}  &d_{i}  \\
             \end{array}
           \right)
,$ $i=1,2. $ \\
  $\tilde{A}_{i}$  &$\tilde{A}_{i} :=\left(
             \begin{array}{cc}
               a_{i} &b_{i} \\
              \frac{ c_{i}}{2} -\frac{1}{2b_{i}} &\frac{ d_{i}}{2}  \\
             \end{array}
           \right)$, $ i=1,2.$\\
\hline\cline{1-2}
\end{tabular}
\end{center}
\label{tab:notation}
\end{table}
\section{Preliminaries}\label{covsec2}

\subsection{Quaternion algebra }
Hamilton invented the 4-D quaternion algebra $\mathbb{H}$ \cite{hamilton} in 1843,
\begin{eqnarray}\label{h0}
\mathbb{H} := \{q=q_r+{\bm \mu}q_\mu+{\bm \nu}q_\nu+{\bm \eta}q_\eta \, | \, q_r,q_\mu,q_\nu,q_\eta\in\mathbb{R}\},
\end{eqnarray}
where $\bm{\mu, \nu, \eta }$ are three units and obey the Hamiliton's multiplication rules:
$
 {\bm \mu}^2= {\bm \nu}^2={\bm \eta}^2=\bm{ \mu\nu\eta}=-1, \bm{ \mu\nu=-\nu\mu=\eta}.
$
The quaternion algebra is a generalization of the complex algebra $\mathbb{C}$.
When $ q_r=0,$ we call $\underline{q}:={\bm \mu}q_\mu+{\bm \nu}q_\nu+{\bm \eta}q_\eta$ a pure imaginary quaternion. Let  $ \mathbb{U}$  be the unit pure imaginary, that is $\mathbb{U}:=\{ q| q=q_r+\underline{q}, q_r=0, \underline{q}^{2}=-1\}.$
To proceeding further, let $\mathbb{H}(\bm{\mu})$ be the field spanned by $ \{ 1, \bm{\mu}\}$, which is the slice quaternion   of $\mathbb{H} $. That is,
 $$\mathbb{H}(\bm{\mu}):=\{ q | q=q=q_{r}+\bm{\mu}q_{\mu}, q_{r}, q_{\mu} \in \mathbb{R}, \bm{\mu} \in \mathbb{U}\}.$$
Meanwhile, the {\it conjugate} of  $q\in \mathbb{H}$ is defined by
$
\overline{q} :=q_r-\left({\bm \mu}q_\mu+{\bm \nu}q_\nu+{\bm \eta}q_\eta\right).\label{h1}
$
The modulus of $q\in\mathbb{H}$ can be defined as
$
|q| := \sqrt{q\overline{q}} = \sqrt{\overline{q}q} = \sqrt{q_r^2+q_\mu^2+q_\nu^2+q_\eta^2}.
$
By $(\ref{h0})$, a quaternionic function $f:\mathbb{R}^2\to\mathbb{H}$ can be expressed in the following form:

\begin{eqnarray}
 f &=&f_{r}+ \bm{\mu}f_{\mu}+\bm{\nu}f_{\nu}+\bm{\eta}f_{\eta},\nonumber\\
   &=&f_{a}+f_{b}\bm{\nu} =f_{a}+\bm{\nu}\overline{f_{b}},
   \label{ab1}
\end{eqnarray}
 where $f_{a},f_{b}\in \mathbb{H(\bm{\mu})} $.



Let $L^{p}(\mathbb{R}^2 , \mathbb{H})$ with integer $p\geq 1$, denote the linear space consisting of all quaternionic functions in $\mathbb{R}^2$ such that
\begin{equation*}
   \Vert f \Vert_{p}:= \int_{\mathbb{R}^2} |f(x,y)|^p dxdy<\infty.
\end{equation*}
That is,
\begin{eqnarray*}
L^{p}(\mathbb{R}^2, \mathbb{H}):=\{ f| f:\mathbb{R}^2 \to \mathbb{H}, \Vert f \Vert_{p}:=\left (\int_{\mathbb{R}^2} |f(x,y)|^p dxdy \right )^{\frac{1}{p}} <\infty \}.
\end{eqnarray*}
If integers $q, p\geq 1,\frac{1}{p}+\frac{1}{q}=1,f\in L^{p}(\mathbb{R}^2, \mathbb{H}) ,g\in L^{q}(\mathbb{R}^2, \mathbb{H})$,
then  H\"{o}lder's   inequality  yields
\begin{eqnarray}\label{hold}
\int_{\mathbb{R}^{2}}|f(x,y)g(x,y)|dxdy \leq \|f(x,y)\|_{p} \|g(x,y)\|_{q}.
\end{eqnarray}
Due to the non-commutativity of quaternion, there are various types of  QFT \cite{ell2014quaternion}, among them the two-sided QFT is defined by

\begin{eqnarray}\label{776}
\mathcal{F}[f](u, v):=\int_{\mathbb{R}^2} e^{-{ \bm\mu}ux}f(x,y)e^{{- \bm\nu}vy}dxdy.\label{hu24}
\end{eqnarray}
where $f\in L^{1}(\mathbb{R}^2, \mathbb{H})$, $(u, v),\in \mathbb{R}^2$ .

\subsection{QLCT}
 In this paper,
we study the  convolution  theorems associated with the two-sided QLCT (hereinafter referred to as the QLCT), which is a generalization of QFT in Eq. (\ref{776}). The  QLCT is defined as following \cite{kou2013uncertainty} , let $f\in L^{1}(\mathbb{R}^2, \mathbb{H})$,
\begin{eqnarray}\label{qlct1}
 \mathcal{L}^{\bm{\mu},\bm{\nu}}_{A_{1},A_{2}}[f](u,v):=  \left\{
\begin{array}{llll}
\int_{\mathbb{R}^2} K_{A_{1}}^{\bm{\mu}}(x,u)f(x,y)K_{A_{2}}^{\bm{\nu}}(y,v)dxdy
& b_{1}, b_{2}\neq 0,\\[1.5ex]
\int_{\mathbb{R}} \sqrt{d_{1}} e^{\bm{\mu}\frac{c_{1}d_{1}u^{2}}{2}}f(du,y)K_{A_{2}}^{\bm{\nu}}(y,v)dy
 & b_{1}=0, b_{2}\neq 0, \\[1.5ex]
\int_{\mathbb{R}} K_{A_{1}}^{\bm{\mu}}(x,u)f(x,dv)\sqrt{d_{2}} e^{\bm{\nu}\frac{c_{2}d_{2}v^{2}}{2}}dx
& b_{1}\neq 0, b_{2} =0, \\[1.5ex]
\sqrt{d_{1}} e^{\bm{\mu}\frac{c_{1}d_{1}u^{2}}{2}}f(du,dv)\sqrt{d_{2}} e^{\bm{\nu}\frac{c_{2}d_{2}v^{2}}{2}}
& b_{1}= 0, b_{2} =0,
\end{array}\right.
\end{eqnarray}
\noindent where  $A_{i}=\left(
             \begin{array}{cc}
               a_{i} &b_{i} \\
               c_{i}  &d_{i}  \\
             \end{array}
           \right)
 \in \mathbb{R}^{2\times 2},$ $i=1,2 $
 both are real matrices parameters with unit determinant, i.e. $det(A_{i})$=$ a_{i}d_{i}-c_{i}b_{i}=1, $ for $i=1,2,$
 and the kernels $K_{A_{1}}^{ \bm{\mu}}(x,u)$ and $K_{A_{2}}^{ \bm{\nu}}(y,v)$ of the QLCT are
\begin{eqnarray}\label{kernel}
K_{A_{1}}^{\bm{\mu}}(x,u):=\frac{1}{ \sqrt{\bm{\mu}2\pi b_{1}}}e^{\bm{\mu}(\frac{a_{1}}{2b_{1}}x^{2}-\frac{1}{b_{1}}ux+\frac{d_{1}}{2b_{1}}u^{2} )}, \quad
K_{A_{2}}^{\bm{\nu}}(y,v):=\frac{1}{\sqrt{\bm{\nu}2\pi b_{2}}}e^{\bm{\nu}(\frac{a_{1}}{2b_{2}}y^{2}-\frac{1}{b_{2}}yv+\frac{d_{2}}{2b_{2}}v^{2} )},
\end{eqnarray}
respectively.  Let $\mathcal{L}^{\bm{\mu},\bm{\nu}}[f]: =\mathcal{L}^{\bm{\mu},\bm{\nu}}_{A_{1},A_{2}}[f]$.
 \par
When  $A_{1}=A_{2}=\left(
             \begin{array}{cc}
               0 &1 \\
               -1  &0  \\
             \end{array}
           \right)$,
the QLCT   reduces to the two-sided  QFT defined in formula (\ref{776}).
When the  kernels $  K_{A_{1}}^{\bm{\mu}}(x,u)$ and  $K_{A_{2}}^{\bm{\nu}}(y,v)$ have the same unit pure imaginary quaternion, the QLCT $\mathcal{L}^{\bm{\mu},\bm{\nu}}_{A_{1},A_{2}}[f](u,v)$ becomes the two-sided QLCT 	$\mathcal{L}^{\bm{\mu},\bm{\mu}}_{A_{1},A_{2}}[f](u,v)$ \cite{bahri2019Two} as follows:
\begin{eqnarray*}
	\mathcal{L}^{\bm{\mu},\bm{\mu}}_{A_{1},A_{2}}[f](u,v):= 
		\int_{\mathbb{R}^2} K_{A_{1}}^{\bm{\mu}}(x,u)f(x,y)K_{A_{2}}^{\bm{\mu}}(y,v)dxdy
\end{eqnarray*}
$\mathcal{L}^{\bm{\mu},\bm{\mu}}_{A_{1},A_{2}}[f](u,v)$ is the special case of the $ \mathcal{L}^{\bm{\mu},\bm{\nu}}_{A_{1},A_{2}}[f](u,v)$. Let
$\mathcal{L}^{\bm{\mu},\bm{\mu}}[f]: =\mathcal{L}^{\bm{\mu},\bm{\mu}}_{A_{1},A_{2}}[f]$.
\par
In this paper, we deal with the case when $b_{1}\neq 0$ and $ b_{2}\neq 0$. If $b_{1}=0$ or $b_{2}= 0,$ $\mathcal{L}^{\bm{\mu},\bm{\nu}}_{A_{1},A_{2}}[f]$ reduces to the special case of one sided QLCT. Moreover, if $b_{1} = b_{2}= 0$,
 $\mathcal{L}^{\bm{\mu},\bm{\nu}}_{A_{1},A_{2}}[f]$  is  just a scaling operation coupled with two chirp multiplications, they are no of interest to our research. 
 \par

 We  conclude this section with two  lemmas, Lemma \ref{L1} tells us that under what conditions the QLCT can be inverted.
\begin{lemma}\cite{hu2016quaternion}\label{L1}
Suppose $f $ and $\mathcal{L}^{\bm{\mu},\bm{\nu}}[f]\in L^{1}(\mathbb{R}^{2}, \mathbb{H})$, then
\begin{eqnarray*}
f(x,y)= \mathcal{L}_{\A_{1}^{-1}, \A_{2}^{-1}}^{\bm{\mu},\bm{\nu}}[\mathcal{L}^{\bm{\mu},\bm{\nu}}[f]](x,y)= \int_{\mathbb{R}^{2}}K_{\A_{1}^{-1}}^{\bm{\mu}}(u,x)\mathcal{L}^{\bm{\mu},\bm{\nu}}[f](u,v)K_{\A_{2}^{-1}}^{\bm{\nu}}(v,y)dudv,
\end{eqnarray*}
for almost everywhere $(x,y)\in \mathbb{R}^2$.
\end{lemma}
The following lemma tells us the  $L^{p}$ bound of the classical convolution operator  $  * $ of quaternionic functions.
Despite its trivial proof, the implications of Lemma \ref{le42} are crucial. As a result, we prove this lemma.
\begin{lemma} \label{le42}
For given functions, $f \in L^{p}(\mathbb{R}^{2}, \mathbb{H}),  g \in L^{1}(\mathbb{R}^{2}, \mathbb{H}),$   integers $ 1 \leq p \leq \infty, $ let
\begin{eqnarray}\label{oldcov1}
(f * g)(x,y):=\int_{R^{2} }f(x-\tau_{1},y-\tau_{2})g(\tau_{1},\tau_{2})d\tau_{1}d\tau_{2},
\end{eqnarray}
then
 \begin{eqnarray}\label{h42}
 \|f * g \|_{p} \leq \|f\|_{p} \|g\|_{1}.
 \end{eqnarray}

\end{lemma}
\begin{nproof}
Firstly, we  need to  show that
\begin{eqnarray*}
\left | \int_{\mathbb{R}^{2}}f(x,y)dxdy \right |\leq \int_{\mathbb{R}^{2}}\left | f(x,y)\right | dxdy,
\end{eqnarray*}
the $ \int_{\mathbb{R}^2} f(x,y)dxdy $ can be rewritten as
\begin{eqnarray*}
\int_{\mathbb{R}^2} f(x,y)dxdy=  \left |\int_{\mathbb{R}^2} f(x,y)dxdy \right |\cdot\Theta,
\end{eqnarray*}
\noindent where $\Theta :=\frac{\int_{\mathbb{R}^2} f(x,y)dxdy}{\left |\int_{\mathbb{R}^2} f(x,y)dxdy \right |}$, if 
$\left |\int_{\mathbb{R}^2} f(x,y)dxdy \right |\neq 0
$,
then
\begin{eqnarray*}
 \int_{\mathbb{R}^2} f(x,y) \Theta^{-1} dxdy=  \left |\int_{\mathbb{R}^2} f(x,y)dxdy \right | \geq 0.
\end{eqnarray*}
Set $w(x,y)=f(x,y) \Theta^{-1} =w_{r}(x,y)+{\bm \mu} w_{\mu}(x,y)+{\bm \nu} w_{\nu}(x,y)+{\bm \eta} w_{\eta}(x,y)$, we have
\begin{eqnarray*}
0 && \leq \left |\int_{\mathbb{R}^2} f(x,y)dxdy \right | =\int_{\mathbb{R}^2} f(x,y) \Theta^{-1} dxdy\\
&&=\int_{\mathbb{R}^2} w_{r}(x,y) dxdy\leq \int_{\mathbb{R}^2}\sqrt{ w^{2}_{r}(x,y)+w^{2}_{\mu}(x,y)+w^{2}_{\nu}(x,y)+w^{2}_{\eta}(x,y)} dxdy\\
 &&= \int_{\mathbb{R}^2}\left | f(x,y) \Theta^{-1}\right| dxdy =\int_{\mathbb{R}^2}\left | f(x,y)\right| dxdy.
 \end{eqnarray*}
Inequality (\ref{h42}) can be obtained by applying the analogous argument to \cite{PWJ2000}.
\end{nproof}

\section{Spatial convolution theorems} \label{covsec3}
 In this section, we introduce  the  spatial  convolution operator  $\circledast^{\bm{\mu},\bm{\nu}}_{A_1,A_2}$. Moreover, the corresponding correlation operator is defined. Finally, some  applications are given at the end of the section.
 
\subsection{Spatial convolution  $\circledast^{{\mu},{\nu}}_{A_1,A_2}  $}
Before we give the definition of the spatial convolution operator  $\circledast^{\bm{\mu},\bm{\nu}}_{A_1,A_2} $.
 Let us  first define  a weight function $W_{A}^{\bm{\mu}}$ from $  \mathbb{R}^2$ to $ \mathbb{H}(\bm{\mu})$ by
 $$ W_{A}^{\bm{\mu}}(t,\tau):=\frac{1}{\sqrt{2\pi b\bm{\mu}}}e^{\bm{\mu}(\tau(\tau-t)\frac{a}{b} )} ,(t,\tau)\in \mathbb{R}^2,  $$ where  $A$ denote the parameter matrix $A=\begin{pmatrix}
 a & b \\
 c & d
\end{pmatrix}$, $a, b,c ,d$ are real numbers satisfying $ad-bc=1.$
%

\begin{definition}(\textbf{Spatial convolution operator $\circledast^{\bm{\mu},\bm{\nu}}_{A_1,A_2} $})\label{def41}
Let $f  \in L^{p}(\mathbb{R}^{2}, \mathbb{H})$,
$g\in L^{1}(\mathbb{R}^{2}, \mathbb{H})$,
integers $ 1 \leq p \leq \infty$,
the spatial  convolution operator  $\circledast^{\bm{\mu},\bm{\nu}}_{A_1,A_2}$ is defined  by
\begin{eqnarray}\label{cov1}
(f\circledast^{\bm{\mu},\bm{\nu}}_{A_1,A_2}  g)(x,y) :=\int_{\mathbb{R }^{2}}W_{A_{1}}^{\bm{\mu}}(x,\tau_{1} )f(\tau _{1},\tau _{2})g(x-\tau _{1},y-\tau _{2})W_{A_{2}}^{\bm{\nu}}(y,\tau _{2})d\tau _{1}d\tau _{2}. \label{hu41}
\end{eqnarray}
\end{definition}
\begin{remark}
When $ a_{i}=d_{i}=0, b_{i}=1$, convolution operator $\circledast^{\bm{\mu},\bm{\nu}}_{A_1,A_2}$ reduces to the classical convolution operator $ *$ defined in Eq. $(\ref{oldcov1})$.
\end{remark}

%

\begin{theorem}($L^{p}\textbf{bounded of spatial convolution}$)\label{th54}
For given functions, $f\in L^{p}(\mathbb{R}^{2}, \mathbb{H})$ and  $g \in L^{1}(\mathbb{R}^{2}, \mathbb{H}),$   integers $ 1 \leq p \leq \infty, $
then
 \begin{eqnarray*}
 \|   f \circledast^{\bm{\mu},\bm{\nu}}_{A_1,A_2}  g   \|_{p} \leq \|f\|_{p} \|g\|_{1}.
 \end{eqnarray*}
\end{theorem}
\begin{nproof}
As $ | W_{A_{1}}^{\bm{\mu}}|=\frac{1}{\sqrt{2\pi b_{1}}}$ and $| W_{A_{2}}^{\bm{\nu}}|=\frac{1}{\sqrt{2\pi b_{2}}} $,  this theorem is easily obtained from Lemma \ref{le42}.
\end{nproof}

\begin{theorem}($\textbf{Spatial convolution theorem  in}$ $\mathbb{H}$)\label{h1}
	Let $f\in L^{p}(\mathbb{R}^{2},\mathbb{H})$, $g\in L^{1}(\mathbb{R}^{2},\mathbb{H})$ and 
$g=g_{1}+\bm{\mu} g_ {2},g_{1},g_{2}\in  \mathbb{H}(\bm{\nu})$, then
	\begin{eqnarray*}
		\mathcal{L}^{\bm{\mu},\bm{\nu}}[f\circledast^{\bm{\mu},\bm{\nu}}_{A_1,A_2} g](u,v)&=&
		e^{\bm{\mu}(-\frac{d_{1}}{2b_{1}}u^{2})}\mathcal{L}^{\bm{\mu},\bm{\nu}}[\mathcal{L}^{\bm{\mu},\bm{\nu}}[f]g_{1}](u,v)	e^{\bm{\nu}(-\frac{d_{2}}{2b_{2}}v^{2})}\nonumber\\
		&&+e^{\bm{\mu}(-\frac{d_{1}}{2b_{1}}u^{2})}\mathcal{L}^{\bm{\mu},\bm{\nu}}[\mathcal{L}^{\bm{\mu},\bm{\nu}}[f\bm{\mu}]g_{2}](u,v)	e^{\bm{\nu}(-\frac{d_{2}}{2b_{2}}v^{2})}.
	\end{eqnarray*}

\begin{nproof}
	By Theorem \ref{th54},
	$\mathcal{L}^{\bm{\mu},\bm{\nu}}[f\circledast_{A_{1},A_{2}} g]$ is well defined. We have
	\begin{eqnarray*}
		&&\mathcal{L}^{\bm{\mu},\bm{\nu}}[f\circledast^{\bm{\mu},\bm{\nu}}_{A_1,A_2} g](u,v)\\
		&=&\int_{R^{2} }K_{A_{1}}^{\bm{\mu}}(x,u)f\circledast^{\bm{\mu},\bm{\nu}}_{A_1,A_2} g(x,y)K_{A_{2}}^{\bm{\nu}}(y,v)dxdy\\
		&=&\int_{R^{2} }\sqrt{\frac{1}{\bm{\mu} b_{1}2\pi}}
		e^{\bm{\mu}(\frac{a_{1}}{2b_{1}}x^2-\frac{1}{b_{1}}xu+\frac{d_{1}}{2b_{1}}u^{2})}f\circledast^{\bm{\mu},\bm{\nu}}_{A_1,A_2}
		  g(x,y)e^{\bm{\nu}(\frac{a_{2}}{2b_{2}}y^2-\frac{1}{b_{2}}yv+\frac{d_{2}}{2b_{2}}v^{2})}\sqrt{\frac{1}{\bm{\nu}b_{2}2\pi}}dxdy\\
		&=&\int_{R^{2} }\sqrt{\frac{1}{\bm{\mu} b_{1}2\pi}}
		e^{\bm{\mu}(\frac{a_{1}}{2b_{1}}x^2-\frac{1}{b_{1}}xu+\frac{d_{1}}{2b_{1}}u^{2})}\bigg(\int_{R^{2} }\sqrt{\frac{1}{\bm{\mu} b_{1}2\pi}}e^{\bm{\mu}(\tau_{1}(\tau_{1}-x)\frac{a_{1}}{b_{1}}) }
		f(\tau _{1},\tau _{2})\\
		&&g(x-\tau _{1},y-\tau _{2})\sqrt{\frac{1}{\bm{\nu} b_{2}2\pi}}e^{\bm{\mu}(\tau_{2}(\tau_{2}-y)\frac{a_{2}}{b_{2}} )}d\tau_{1}d\tau_{2}\bigg)
		e^{\bm{\nu}(\frac{a_{2}}{2b_{2}}y^2-\frac{1}{b_{2}}yv+\frac{d_{2}}{2b_{2}}v^{2})}\sqrt{\frac{1}{\bm{\nu} b_{2}2\pi}}dxdy.
	\end{eqnarray*}
	By changing variables $s=x-\tau_{1}, t=y-\tau_{2},$ we get
	\begin{eqnarray*}
		&&\mathcal{L}^{\bm{\mu},\bm{\nu}}[f\circledast^{\bm{\mu},\bm{\nu}}_{A_1,A_2}  g](u,v)\\
		&=&\int_{R^{2} }\sqrt{\frac{1}{\bm{\mu}b_{1}2\pi}}
		e^{\bm{\mu}(\frac{a_{1}}{2b_{1}}(s+\tau_{1})^2-\frac{1}{b_{1}}u(s+\tau_{1})+\frac{d_{1}}{2b_{1}}u^{2})}\\
		&& \bigg(\int_{R^{2} }\sqrt{\frac{1}{\bm{\mu}b_{1}2\pi}}e^{\bm{\mu}(\tau_{1}(-s)\frac{a_{1}}{b_{1}}) }f(\tau _{1},\tau _{2})
		(g_1(s,t)+\bm{\mu}g_2(s,t))\sqrt{\frac{1}{\bm{\nu}b_{2}2\pi}}e^{\bm{\nu}(\tau_{2}(-t)\frac{a_{2}}{b_{2}} )}d\tau _{1}d\tau _{2}\bigg)\\
		&&\sqrt{\frac{1}{\bm{\nu} b_{2}2\pi}}e^{\bm{\nu}(\frac{a_{2}}{2b_{2}}(t+\tau_{2})^2-\frac{1}{b_{2}}(t+\tau_{2})v+\frac{d_{2}}{2b_{2}}v^{2})}dsdt\\
		&=&e^{\bm{\mu}(-\frac{d_{1}}{2b_{1}}u^{2})}\mathcal{L}^{\bm{\mu},\bm{\nu}}[\mathcal{L}^{\bm{\mu},\bm{\nu}}[f]g_{1}](u,v)	e^{\bm{\nu}(-\frac{d_{2}}{2b_{2}}v^{2})}\\
		&&+e^{\bm{\mu}(-\frac{d_{1}}{2b_{1}}u^{2})}\mathcal{L}^{\bm{\mu},\bm{\nu}}[\mathcal{L}^{\bm{\mu},\bm{\nu}}[f\bm{\mu}]g_{2}](u,v)	e^{\bm{\nu}(-\frac{d_{2}}{2b_{2}}v^{2})}
	\end{eqnarray*}
	which completes the proof.
\end{nproof}
\end{theorem}
For the purpose of obtaining the applications of the correlation theorem, the  Fredholm integral
equations involving special kernels, some special system of second-order partial differential equations, and  designing  multiplication filters in the QLCT domain, we consider the special case of the spatial convolution theorem with $\mathcal{L}^{\bm{\mu},\bm{\mu}}$.
For simplifying the notations, denote
$$\circledast_{A_1,A_2}:=\circledast^{\bm{\mu},\bm{\mu}}_{A_1,A_2}.$$ We divide the subsections into the following: we first study the spatial convolution in $\mathbb{H}(\bf{\mu})$, then apply it to deduce the result in the general quaternion space $\mathbb{H}$. Firstly, we extend
 the 1D convolution theorem of LCT  \cite{wei2012new} to  $\mathbb{H(\bm{\mu})}$. The following result is different with those derived in \cite{deng2006convolution,wei2009convolution,deng2010comments,wei2012new,wei2012convolution,shi2014generalized}, ours preserves the classical Fourier transform property. 
\subsubsection {$\circledast_{A_1,A_2} $in $\mathbb{H}(\bf{\mu})$} 
Let
$\tilde{A}_{i} :=\left(
\begin{array}{cc}
	a_{i} &b_{i} \\
	\frac{ c_{i}}{2} -\frac{1}{2b_{i}} &\frac{ d_{i}}{2}  \\
\end{array}
\right)$, $ i=1,2.$
It is easy to see that
$$   \mathcal{L}_{\tilde{A}_{1},\tilde{A}_{2} }^{\bm{\mu},\bm{\mu}}[f](u,v  )= e^{-\bm{\mu}(\frac{d_{1}}{4b_{1}}u^{2} )} \mathcal{L} ^{\bm{\mu},\bm{\mu}}[f](u,v)e^{-\bm{\mu}(\frac{d_{2}}{4b_{2}}v^{2} )}.$$

\begin{theorem}(\textbf{Spatial convolution theorem in} $\mathbb{H}(\bm{\mu}))$ \label{le43}
 Let $f\in L^{p}(\mathbb{R}^{2},\mathbb{H}(\bm{\mu}))$ and  $g \in L^{1}(\mathbb{R}^{2}, \mathbb{H}(\bm{\mu})),$
then
\begin{eqnarray}
\mathcal{L}^{\bm{\mu},\bm{\mu}}[f\circledast_{A_{1},A_{2}} g](u,v)&&=
e^{\bm{\mu}(-\frac{d_{1}}{2b_{1}}u^{2})}e^{\bm{\mu}(-\frac{d_{2}}{2b_{2}}v^{2})}\mathcal{L}^{\bm{\mu},\bm{\mu}}[f](u,v)
\mathcal{L}^{\bm{\mu},\bm{\mu}}[g](u,v) \mbox{ and }\label{hu81}\\
\mathcal{L}^{\bm{\mu},\bm{\mu}}[f\circledast_{A_{1},A_{2}} g](u,v)&&=
\mathcal{L}_{\tilde{A}_{1},\tilde{A}_{2} }^{\bm{\mu},\bm{\mu}}[f](u,v )
 \mathcal{L}_{\tilde{A}_{1},\tilde{A}_{2} }^{\bm{\mu},\bm{\mu}}[g]
(u, v).\label{h81}
\end{eqnarray}
\end{theorem}

\begin{nproof}
By Theorem \ref{th54},
$\mathcal{L}^{\bm{\mu},\bm{\mu}}[f\circledast_{A_{1},A_{2}} g](u,v)$ is well defined. We have that
\begin{eqnarray*}
&&\mathcal{L}^{\bm{\mu},\bm{\mu}}[f\circledast_{A_{1},A_{2}} g](u,v)\\
&=&\int_{R^{2} }K_{A_{1}}^{\bm{\mu}}(x,u)f\circledast_{A_{1},A_{2}} g(x,y)K_{A_{2}}^{\bm{\mu}}(y,v)dxdy\\
&=&\int_{R^{2} }\sqrt{\frac{1}{\bm{\mu} b_{1}2\pi}}\sqrt{\frac{1}{\bm{\mu}b_{2}2\pi}}
e^{\bm{\mu}(\frac{a_{1}}{2b_{1}}x^2-\frac{1}{b_{1}}xu+\frac{d_{1}}{2b_{1}}u^{2})}f\circledast_{A_{1},A_{2}} g(x,y)e^{\bm{\mu}(\frac{a_{2}}{2b_{2}}y^2-\frac{1}{b_{2}}yv+\frac{d_{2}}{2b_{2}}v^{2})}dxdy\\
&=&\int_{R^{2} }\sqrt{\frac{1}{\bm{\mu} b_{1}2\pi}}\sqrt{\frac{1}{\bm{\mu} b_{2}2\pi}}
e^{\bm{\mu}(\frac{a_{1}}{2b_{1}}x^2-\frac{1}{b_{1}}xu+\frac{d_{1}}{2b_{1}}u^{2})}\bigg(\int_{R^{2} }\sqrt{\frac{1}{\bm{\mu} b_{1}2\pi}}e^{\bm{\mu}(\tau_{1}(\tau_{1}-x)\frac{a_{1}}{b_{1}}) }
f(\tau _{1},\tau _{2})\\
&& g(x-\tau _{1},y-\tau _{2})\sqrt{\frac{1}{\bm{\mu} b_{2}2\pi}}e^{\bm{\mu}(\tau_{2}(\tau_{2}-y)\frac{a_{2}}{b_{2}} )}d\tau_{1}d\tau_{2}\bigg)
e^{\bm{\mu}(\frac{a_{2}}{2b_{2}}y^2-\frac{1}{b_{2}}yv+\frac{d_{2}}{2b_{2}}v^{2})}dxdy.
\end{eqnarray*}
By changing variables $s=x-\tau_{1}, t=y-\tau_{2},$ we have
\begin{eqnarray*}
&&\mathcal{L}^{\bm{\mu},\bm{\mu}}[f\circledast_{A_{1},A_{2}} g](u,v)\\
&=&\int_{R^{2} }\sqrt{\frac{1}{\bm{\mu}b_{1}2\pi}}\sqrt{\frac{1}{\bm{\mu} b_{2}2\pi}}
e^{\bm{\mu}(\frac{a_{1}}{2b_{1}}(s+\tau_{1})^2-\frac{1}{b_{1}}u(s+\tau_{1})+\frac{d_{1}}{2b_{1}}u^{2})}\\
&&  \bigg(\int_{R^{2} }\sqrt{\frac{1}{\bm{\mu}b_{1}2\pi}}e^{\bm{\mu}(\tau_{1}(-s)\frac{a_{1}}{b_{1}}) }f(\tau _{1},\tau _{2})
 g(s,t)\sqrt{\frac{1}{\bm{\mu}b_{2}2\pi}}e^{\bm{\mu}(\tau_{2}(-t)\frac{a_{2}}{b_{2}} )}d\tau _{1}d\tau _{2}\bigg)\\
&&  e^{\bm{\mu}(\frac{a_{2}}{2b_{2}}(t+\tau_{2})^2-\frac{1}{b_{2}}(t+\tau_{2})v+\frac{d_{2}}{2b_{2}}v^{2})}dsdt\\
&=&e^{-\bm{\mu}(\frac{d_{1}}{2b_{1}}u^{2}+\frac{d_{2}}{2b_{2}}v^{2})}\mathcal{L}^{(\bm{\mu},\bm{\mu})}[f](u,v)\mathcal{L}^{(\bm{\mu},\bm{\mu})}[g](u,v)\\
&=&  \mathcal{L}_{\tilde{A}_{1},\tilde{A}_{2}}^{\bm{\mu},\bm{\mu}}[f](u,v  ) \mathcal{L}_{\tilde{A}_{1},\tilde{A}_{2}}^{\bm{\mu},\bm{\mu}}[g](u,v  ).
\end{eqnarray*}
It completes the proof.
\end{nproof}

\begin{remark}
\begin{itemize}
  \item Eqs. (\ref{hu81}) and (\ref{h81}) are identitial in  mathematics. The reader can refer to Figures  \ref{ff1} and \ref{ff2} for details on how they can be used to design different multiplicative filters.
  
 \item Eq.  (\ref{h81}) of Theorem \ref{le43} generalizes the 1D convolution theorem of the LCT in paper \cite{wei2012new} to quaternion setting $\mathbb{H}(\bm{\mu})$. Eq. (\ref{h81}) shows that the  QLCT of  convolution of two  quaternionic functions is the product of their QLCTs without the chip signals. If $A_1=A_2=\left(  \begin{array}{cc}
 	1 &0 \\
 	0  &1  \\
 \end{array} \right) $, then Eq. (\ref{h81}) reduces to the classical case of the 2-D Fourier transform.

\end{itemize}
\end{remark}
The  $\textbf{ Parseval's  formula}$ and $ \textbf{Energy  theorem}$ in terms of convolution operator $\circledast_{A_1,A_2}$ are given as follows.
\begin{theorem}(\textbf{Parseval's formula})\label{le410}
Let $f,g,$ $\mathcal{L}^{\bm{\mu},\bm{\mu}}[f]$ and $\mathcal{L}^{\bm{\mu},\bm{\mu}}[g]\in L^{1}\bigcap L^{2}(\mathbb{R}^2, \mathbb{H}(\bm{\mu})) $, then
 \begin{eqnarray*}
 \int_{\mathbb{R}^{2}}f(\tau_{1},\tau_{2})\overline{ r(\tau_{1},\tau_{2})}d\tau_{1}d\tau_{2}= \int_{\mathbb{R}^{2}}\mathcal{L}^{\bm{\mu},\bm{\mu}}[f](u,v)
 \overline{\mathcal{L}^{\bm{\mu},\bm{\mu}}[r](u,v)}dudv,
 \end{eqnarray*}
 \noindent where $\overline{r( \tau_{1},\tau_{2})}=e^{\bm{\mu}\tau_{1}^{2}\frac{a_{1}}{b_{1}}}g(-\tau_{1},-\tau_{2})e^{\bm{\mu}\tau_{2}^{2}\frac{a_{2}}{b_{2}}}$.\end{theorem}
\begin{nproof} According to Theorem \ref{th54},
we see that $f\circledast_{A_1,A_2} g \in L^{1}(\mathbb{R}^2, \mathbb{H}(\bm{\mu}))$,  then $ \mathcal{L}^{\bm{\mu},\bm{\mu}}[f\circledast_{A_1,A_2} g]$ is well defined.  Moreover,    $ \mathcal{L}^{\bm{\mu},\bm{\mu}}[f\circledast_{A_1,A_2} g]\in L^{1}(\mathbb{R}^2, \mathbb{H}(\bm{\mu})) $ is followed from the H\"{o}lder's    inequality $(\ref{hold})$,  and
 $\mathcal{L}^{\bm{\mu} \bm{\mu}}[f],$ $\mathcal{L}^{\bm{\mu},\bm{\mu}}[g]$ $\in L^{1}\bigcap  L^{2}(\mathbb{R}^2, \mathbb{H}(\bm{\mu})).$
Due to Lemma \ref{L1}, we can recover the $(f\circledast_{A_1,A_2} g)$ from its QLCT. That is, 
\begin{eqnarray*}
&&\int_{\mathbb{R}^{2} }W_{A_{1}}^{\bm{\mu}}(x,\tau_{1} )f(\tau _{1},\tau _{2})g(x-\tau _{1},y-\tau _{2})W_{A_{2}}^{\bm{\mu}}(y,\tau _{2})d\tau _{1}d\tau _{2}\\
=&&\int_{\mathbb{R}^{2}}\frac{1}{\sqrt{-\bm{\mu}2\pi b_{1}}}e^{\bm{\mu}( \frac{-d_{1}}{2b_{1}}u^{2} +\frac{1}{b_{1}}ux-  \frac{-a_{1}}{2b_{1}}x^{2}  )}
e^{-\bm{\mu}\frac{d_{1}}{2b_{1}}u^{2}}e^{-\bm{\mu}\frac{d_{2}}{2b_{2}}v^{2}}
\mathcal{L}^{\bm{\mu},\bm{\mu}}[f](u,v)\mathcal{L}^{\bm{\mu},\bm{\mu}}[g](u,v)\\
&& \frac{1}{\sqrt{-\bm{\mu}2\pi b_{2}}}e^{\bm{\mu}( \frac{-d_{2}}{2b_{2}}v^{2} +\frac{1}{b_{2}}vy-  \frac{-a_{2}}{2b_{2}}y^{2}  )}
dudv.
\end{eqnarray*}
Setting $x=y=0 $ in the above equation, we obtain
\begin{eqnarray*}
&&\int_{\mathbb{R}^{2} }e^{\bm{\mu}\tau_{1}^{2}\frac{a_{1}}{b_{1}}}f(\tau _{1},\tau _{2})g(-\tau _{1},-\tau _{2})e^{\bm{\mu}\tau_{2}^{2}\frac{a_{2}}{b_{2}}}d\tau _{1}d\tau _{2}\\
=&&\int_{\mathbb{R}^{2}}e^{\bm{\mu} \frac{-d_{1}}{b_{1}}u^{2}  }
\mathcal{L}^{\bm{\mu},\bm{\mu}}[f](u,v)\mathcal{L}^{\bm{\mu},\bm{\mu}}[g](u,v)
e^{\bm{\mu} \frac{-d_{2}}{b_{2}}v^{2}  }
dudv.
\end{eqnarray*}
Note that $\overline{r(\tau_{1},\tau_{2})}=e^{\bm{\mu}\tau_{1}^{2}\frac{a_{1}}{b_{1}}}g(-\tau_{1},-\tau_{2})e^{\bm{\mu}\tau_{2}^{2}\frac{a_{2}}{b_{2}}}$,
by straightforward computation, we have
\begin{eqnarray*}
&&e^{\bm{\mu} \frac{-d_{1}}{b_{1}}u^{2}  }\mathcal{L}^{\bm{\mu},\bm{\mu}}[g](u,v)e^{\bm{\mu} \frac{-d_{2}}{b_{2}}v^{2}  }\\
=&&e^{\bm{\mu} \frac{-d_{1}}{b_{1}}u^{2}  }\int_{\mathbb{R}^{2}} \frac{1}{\sqrt{\bm{\mu} 2\pi b_{1}}}
e^{\bm{\mu}(\frac{a_{1}}{2b_{1}}\tau_{1}^2-\frac{1}{b_{1}}\tau_{1}u+\frac{d_{1}}{2b_{1}}u^{2})}
e^{-\bm{\mu}\tau_{1}^{2}\frac{a_{1}}{b_{1}}}\overline{r}(-\tau_{1},-\tau_{2})e^{-\bm{\mu}\tau_{2}^{2}\frac{a_{2}}{b_{2}}}\\
&& \frac{1}{\sqrt{\bm{\mu}2\pi b_{2}}}  e^{\bm{\mu}(\frac{a_{2}}{2b_{2}}\tau_{2}^2-\frac{1}{b_{2}}\tau_{2}v+\frac{d_{2}}{2b_{2}}v^{2})}d\tau_{1}d\tau_{2}
e^{\bm{\mu} \frac{-d_{2}}{b_{2}}v^{2}  }\\
=&&\int_{\mathbb{R}^{2}} \frac{1}{\sqrt{\bm{\mu}2\pi b_{1}}}
e^{\bm{\mu}(-\frac{a_{1}}{2b_{1}}\tau_{1}^2+\frac{1}{b_{1}}\tau_{1}u-\frac{d_{1}}{2b_{1}}u^{2})}
\overline{r}(\tau_{1},\tau_{2})
\frac{1}{\sqrt{\bm{\mu}2\pi b_{2}}}  e^{\bm{\mu}(-\frac{a_{2}}{2b_{2}}\tau_{2}^2+\frac{1}{b_{2}}\tau_{2}v-\frac{d_{2}}{2b_{2}}v^{2})}d\tau_{1}d\tau_{2}\\
=&&\int_{\mathbb{R}^{2}} \frac{1}{\sqrt{-\bm{\mu}2\pi b_{1}}}
e^{\bm{\mu}(-\frac{a_{1}}{2b_{1}}\tau_{1}^2+\frac{1}{b_{1}}\tau_{1}u-\frac{d_{1}}{2b_{1}}u^{2})}
\overline{r}(\tau_{1},\tau_{2})
\frac{1}{\sqrt{-\bm{\mu}2\pi b_{2}}}  e^{\bm{\mu}(-\frac{a_{2}}{2b_{2}}\tau_{2}^2+\frac{1}{b_{2}}\tau_{2}v-\frac{d_{2}}{2b_{2}}v^{2})}d\tau_{1}d\tau_{2}\\
=&&\overline{\mathcal{L}^{\bm{\mu},\bm{\mu}}[r](u,v)},
\end{eqnarray*}
which completes the proof.
\end{nproof}

Moreover, if $r(x,y)=f(x,y)$, then we obtain the following $ \textbf{Energy  theorem}$.
\begin{corollary} (\textbf{Energy  theorem})\label{te45}
 Suppose  that  $f $ and $\mathcal{L}^{\bm{\mu},\bm{\mu}}[f] \in L^{1}\bigcap  L^{2}(\mathbb{R}^2, \mathbb{H}(\bm{\mu})),$ then
 \begin{eqnarray*}
 \int_{\mathbb{R}^{2}}\left |f(\tau_{1},\tau_{2})\right |^{2} d\tau_{1}d\tau_{2}= \int_{\mathbb{R}^{2}}\left |\mathcal{L}^{\bm{\mu},\bm{\mu}}[f](u,v)\right |^{2}dudv.
 \end{eqnarray*}

\end{corollary}


\begin{corollary}
 (\textbf{Product  theorem on $\mathbb{H}(\bm{\mu})$}) \label{le44}
If $f,g \in L^{1}\bigcap L^{2}(\mathbb{R}^2,\mathbb{H}(\bm{\mu})),$ and
$\mathcal{L}^{\bm{\mu},\bm{\mu}}[f]\in L^{p}(\mathbb{R}^2,\mathbb{H}(\bm{\mu})),\mathcal{L}^{\bm{\mu},\bm{\mu}}[g]\in L^{1}(\mathbb{R}^2,\mathbb{H}(\bm{\mu})),$
then
\begin{eqnarray}
\mathcal{L}^{\bm{\mu},\bm{\mu}}{[ e^{\bm{\mu}\frac{a_{1}}{2b_{1}}x^{2}}e^{\bm{\mu}\frac{a_{2}}{2b_{2}}y^{2}} f g]}(u,v)&&=\mathcal{L}^{\bm{\mu},\bm{\mu}}[f]
\circledast_{A_{1}^{-1},A_{2}^{-1}}  \mathcal{L}^{\bm{\mu},\bm{\mu}}[g](u,v), \label{hu82}\\
\mathcal{L}_{\breve{A}_{1},  \breve{A}_{2}}^{\bm{\mu},\bm{\mu}}[ f g](u,v)&&=\mathcal{L}^{\bm{\mu},\bm{\mu}}[f]
\circledast_{A_{1}^{-1},A_{2}^{-1}} \mathcal{L}^{\bm{\mu},\bm{\mu}}[g](u,v), \label{h82}
\end{eqnarray}
where  $\breve{A}_{i} :=\left(
             \begin{array}{cc}
               2a_{i} &b_{i} \\
               \frac{1}{b_{i}} +2c_{i} & d_{i}  \\
             \end{array}
           \right),  i=1,2.$
\end{corollary}
\begin{nproof}
By
$$\int _{\mathbb{R}^2}\mid f(x,y)g(x,y)\mid dxdy \leq \parallel f(x,y) \parallel^{2} _{2}\parallel g(x,y) \parallel^{2} _{2},$$
we have
$\mathcal{L}^{\bm{\mu},\bm{\mu}}
[ e^{\bm{\mu}\frac{a_{1}}{2b_{1}}x^{2}}e^{\bm{\mu}\frac{a_{2}}{2b_{2}}y^{2}} fg]
(u,v)$ is well defined.
Let $L(u,v):=\mathcal{L}^{\bm{\mu},\bm{\mu}}[f] \circledast_{A_{1}^{-1},A_{2}^{-1}}\mathcal{L}^{\bm{\mu},\bm{\mu}}[g](u,v),$
it follows from Theorem \ref{th54} that $L \in L^{1}(\mathbb{R}^2,\mathbb{H}(\bm{\mu}))$. Denote
\begin{eqnarray*}
I &:=&\int_{\mathbb{R}^2 }K_{A^{-1}_{1}}^{\bm{\mu}}(u,x) L(u,v)K_{A^{-1}_{2}}^{\bm{\mu}}(v,y)dudv\\
&=&\int _{\mathbb{R}^2}\frac{1}{\sqrt{-2\bm{\mu}b_{1}\pi }}e^{\bm{\mu}(-\frac{d_{1}}{2b_{1}}u^2+\frac{1}{b_{1}}ux-\frac{a_{1}}{2b_{1}}x^{2})}\\
&& \int _{\mathbb{R}^2}\frac{1}{\sqrt{-2\bm{\mu}b_{1}\pi }}e^{\bm{\mu}\tau _{1}(\tau _{1}-u)(-\frac{d_{1}}{b_{1}})}
\mathcal{L}^{\bm{\mu},\bm{\mu}}[f](u-\tau_{1} ,v-\tau_{2})\mathcal{L}^{\bm{\mu},\bm{\mu}}[g](\tau_{1}, \tau_{2})\frac{1}{\sqrt{-2\bm{\mu} b_{2}\pi }}\\&&  e^{\bm{\mu}\tau _{2}(\tau _{2}-v)(-\frac{d_{2}}{b_{2}})}d\tau_{1} d\tau_{2}
\frac{1}{\sqrt{-2\bm{\mu} b_{2}\pi }}e^{\bm{\mu}(-\frac{d_{2}}{2b_{2}}v^2+\frac{1}{b_{2}}vy-\frac{a_{2}}{2b_{2}}y^{2})}dudv,
\end{eqnarray*}
 substituting $u-\tau_{1}=s, v-\tau_{2}=t$ into the above equation, we have
 \begin{eqnarray*}
 I&=&\int _{\mathbb{R}^2}\frac{1}{\sqrt{-2\bm{\mu} b_{1}\pi }}\frac{1}{\sqrt{-2\bm{\mu} b_{2}\pi }}e^{\bm{\mu}(-\frac{d_{1}}{2b_{1}}(s+\tau_{1})^2+\frac{1}{b_{1}}(s+\tau_{1})x-\frac{a_{1}}{2b_{1}}x^{2})}\\
&&\int _{\mathbb{R}^2}\frac{1}{\sqrt{-2\bm{\mu} b_{1}\pi }}e^{\bm{\mu} \tau _{1}s\frac{d_{1}}{b_{1}}}
\mathcal{L}^{\bm{\mu},\bm{\mu}}[f](s,t)\mathcal{L}^{\bm{\mu},\bm{\mu}}[g](\tau_{1} ,\tau_{2})\frac{1}{\sqrt{-2\bm{\mu} b_{2}\pi }}e^{\bm{\mu}\tau _{2}s\frac{d_{2}}{b_{2}}}d\tau_{1} d\tau_{2}\\
&& e^{\bm{\mu}(-\frac{d_{2}}{2b_{2}}(t+\tau_2)^2+\frac{1}{b_{2}}y(t+\tau_{2})-\frac{a_{2}}{2b_{2}}y^{2})}dsdt\\
&=&e^{\bm{\mu}\frac{a_{1}}{2b_{1}}x^{2}}e^{\bm{\mu}\frac{a_{2}}{2b_{2}}y^{2}} f(x,y)g(x,y).
 \end{eqnarray*}
This gives the proof of Eq. (\ref{hu82}).
  The proof of the Eq.$ (\ref{h82})$ follows by a similar argument,
\begin{eqnarray*}
&&\int_{\mathbb{R}^2 }K_{\breve{A}^{-1}_{1}}^{\bm{\mu}}(u,x) L(u,v)K_{\breve{A}^{-1}_{2}}^{\bm{\mu}}(v,y)dudv\\
&=&\int _{\mathbb{R}^2}\frac{1}{\sqrt{-2\bm{\mu}b_{1}\pi }}e^{\bm{\mu}(-\frac{d_{1}}{2b_{1}}u^2+\frac{1}{b_{1}}ux-\frac{a_{1}}{b_{1}}x^{2})}
\int _{\mathbb{R}^2}\frac{1}{\sqrt{-2\bm{\mu}b_{1}\pi }}e^{\bm{\mu}\tau _{1}(\tau _{1}-u)(-\frac{d_{1}}{b_{1}})}
\mathcal{L}^{\bm{\mu},\bm{\mu}}[f](u-\tau_{1} ,v-\tau_{2})\\
&&\mathcal{L}^{\bm{\mu},\bm{\mu}}[g](\tau_{1}, \tau_{2})\frac{1}{\sqrt{-2\bm{\mu} b_{2}\pi }}e^{\bm{\mu}\tau _{2}(\tau _{2}-v)(-\frac{d_{2}}{b_{2}})}d\tau_{1} d\tau_{2}
\frac{1}{\sqrt{-2\bm{\mu} b_{2}\pi }}e^{\bm{\mu}(-\frac{d_{2}}{2b_{2}}v^2+\frac{1}{b_{2}}vy-\frac{a_{2}}{b_{2}}y^{2})}dudv\\
&=&f(x,y)g(x,y).
\end{eqnarray*}
 It completes the proof.

\end{nproof}
\begin{remark} \label{re45}
\begin{itemize}
  \item If the functions take value in $\mathbb{H}(\bm{\mu})$,
from Eq. (\ref{h81}), the convolution operation  in the spatial domain becomes  the product operation when applying  QLCT, and vice versa. According to Eq.(\ref{h82}), the dual convolution operation in the QLCT domain is converted to  the product operation in the spatial domain  when applying inverse  QLCT.
Hence Eqs.  (\ref{h82})  and   (\ref{h81})   are particularly useful in filter design in the spatial domain and the  QLCT domain, respectively.
\item
 If      $A_{1}=A_{2} =\begin{pmatrix}
 0&1 \\
 -1& 0
\end{pmatrix}$, then the convolution Theorems \ref{le43} and Corollary \ref{le44}   reduce to the classical case in 2-D Fourier transform domain as follows:
\begin{eqnarray*}
&&\int _{\mathbb{R}^2}f(\tau _{1},\tau _{2})g(x-\tau _{1},y-\tau _{2})d\tau _{1}d\tau _{2}\Leftrightarrow \hat{f}(u,v)\hat{g}(u,v),\\
&&\int _{\mathbb{R}^2}\hat{f}(\tau _{1},\tau _{2})\hat{g}(u-\tau _{1},v-\tau _{2})d\tau _{1}d\tau _{2}\Leftrightarrow f(x,y)g(x,y),
\end{eqnarray*}
\noindent where $\hat{f}$ and $\hat{g}$ are 2-D Fourier transform of $f$ and $ g$, respectively.
\end{itemize}
\end{remark}

\subsubsection{$\circledast_{A_1,A_2} $in  $\mathbb{H}$}
In this subsection, we consider the  spatial convolution theorem in $ \mathbb{H}.$
\begin{theorem}(\textbf{Spatial convolution theorem in $\mathbb{H}$}) \label{th46}
For two given quaternionic functions, $f\in L^{p}(\mathbb{R}^{2},\mathbb{H}), g\in L^{1}(\mathbb{R}^{2},\mathbb{H}),$
$f=f_{a}+f_{b}\bm{\nu},g=g_{a}+\bm{\nu}\overline{ g_{b}},$ with $ f_{a},f_{b},g_{a},g_{b} \in \mathbb{H}(\bm{\mu})$ which is defined in Eq.(\ref{ab1}),
then
\begin{eqnarray}
&&\mathcal{L}^{\bm{\mu},\bm{\mu}}[(f
\circledast_{A_{1},A_{2}} g)](u,v) \nonumber\\
=&&e^{-\bm{\mu}(\frac{d_{1}}{2b_{1}}u^2+\frac{d_{2}}{2b_{2}}v^2)}
\big[\mathcal{L}^{\bm{\mu},\bm{\mu}}[f_{a}](u,v)\mathcal{L}^{\bm{\mu},\bm{\mu}}[g_{a}](u,v)
-\mathcal{L}^{\bm{\mu},\bm{\mu}}[f_{b}](u,v)\mathcal{L}^{\bm{\mu},\bm{\mu}}[\overline{ g_{b}}](u,v)\big]\nonumber\\
&&+e^{-\bm{\mu}(\frac{d_{1}}{2b_{1}}u^2-\frac{d_{2}}{2b_{2}}v^2)}
\big[\mathcal{L}^{\bm{\mu},-\bm{\mu}}[f_{b}](u,v)\mathcal{L}^{\bm{\mu},-\bm{\mu}}[\overline{g}_{a}](u,v)
+\mathcal{L}^{\bm{\mu},-\bm{\mu}}[f_{a}](u,v)\mathcal{L}^{\bm{\mu},-\bm{\mu}}[g_{b}](u,v)\big]\bm{\nu}\label{hu83}
\end{eqnarray}
and
\begin{eqnarray}
&&\mathcal{L}^{\bm{\mu},\bm{\mu}}[(f\circledast_{A_{1},A_{2}}  g)](u,v)\nonumber\\
=&&
[\mathcal{L}_{\tilde{A}_{1},\tilde{A}_{2}}^{\bm{\mu},\bm{\mu}}[f_{a}](u,v )\mathcal{L}_{\tilde{A}_{1},\tilde{A}_{2}}^{\bm{\mu},\bm{\mu}}[g_{a}](u,v)
-\mathcal{L}_{\tilde{A}_{1},\tilde{A}_{2}}^{\bm{\mu},\bm{\mu}}[f_{b}](u,v)
\mathcal{L}_{\tilde{A}_{1},\tilde{A}_{2}}^{\bm{\mu},\bm{\mu}}[\overline{ g_{b}}](u,v)] \nonumber\\
&&+
\big[\mathcal{L}^{\bm{\mu},-\bm{\mu}}_{\tilde{A}_{1},\tilde{A}_{2}}[f_{b}](u,v)\mathcal{L}^{\bm{\mu},-\bm{\mu}}_{\tilde{A}_{1},\tilde{A}_{2}}[\overline{g}_{a}](u,v)
+\mathcal{L}^{\bm{\mu},-\bm{\mu}}_{\tilde{A}_{1},\tilde{A}_{2}}[f_{a}](u,v)\mathcal{L}^{\bm{\mu},-\bm{\mu}}_{\tilde{A}_{1},\tilde{A}_{2}}[g_{b}](u,v)\big]\bm{\nu}. \label{h83}
\end{eqnarray}
\end{theorem}

\begin{nproof}
By Theorem \ref{th54}, $(f\circledast_{A_{1},A_{2}}  g)(x,y)$ is well defined and belongs to $L^{1}(\mathbb{R}^{2},\mathbb{H}).$\\
Since
\begin{eqnarray*}
f\circledast g(x,y)
=&&(f_{a}+f_{b}\bm{\nu})\circledast_{A_{1},A_{2}}  (g_{a}+\bm{\nu}\overline{ g_{b}})(x,y),\\
=&&f_{a}\circledast_{A_{1},A_{2}}  g_{a}(x,y )-f_{b}\circledast_{A_{1},A_{2}}  \overline{ g_{b}}(x,y )+f_{b}\circledast_{A_{1},A_{2}}  \bm{\nu} g_{a}(x,y )+f_{a}\circledast_{A_{1},A_{2}}  \bm{\nu} \overline{ g_{b}}(x,y ),
\end{eqnarray*}
then
\begin{eqnarray}\label{huc1}
&&(f\circledast_{A_{1},A_{2}}  g)(x,y)\nonumber\\
&=&\int_{-\infty }^{\infty}\int_{-\infty }^{\infty}W_{A_{1}}^{\bm{\mu}}(x,\tau_{1})f(\tau_{1},\tau_{2})g(x-\tau_{1},y-\tau_{2})W_{A_{2}}^{\bm{\mu}}(y,\tau_{2})d\tau_{1}d\tau_{2} \nonumber\\
&=&\bigg[  \int_{-\infty }^{\infty}\int_{-\infty }^{\infty}W_{A_{1}}^{\bm{\mu}}(x,\tau_{1})f_{a}(\tau_{1},\tau_{2})g_{a}(x-\tau_{1},y-\tau_{2})W_{A_{2}}^{\bm{\mu}}(y,\tau_{2})d\tau_{1}d\tau_{2} \\
&&-\int_{-\infty }^{\infty}\int_{-\infty }^{\infty}W_{A_{1}}^{\bm{\mu}}(x,\tau_{1})f_{b}(\tau_{1},\tau_{2})\overline{g_{b}}(x-\tau_{1},y-\tau_{2})W_{A_{2}}^{\bm{\mu}}(y,\tau_{2})d\tau_{1}d\tau_{2}\bigg ] \nonumber \\
&&+\bigg[ \int_{-\infty }^{\infty}\int_{-\infty }^{\infty}W_{A_{1}}^{\bm{\mu}}(x,\tau_{1})f_{b}(\tau_{1},\tau_{2})\overline{g_{a}}(x-\tau_{1},y-\tau_{2})\overline{W_{A_{2}}^{\bm{\mu}}}(y,\tau_{2})d\tau_{1}d\tau_{2} \nonumber\\
&&+\int_{-\infty }^{\infty}\int_{-\infty }^{\infty}W_{A_{1}}^{\bm{\mu}}(x,\tau_{1})f_{a}(\tau_{1},\tau_{2})g_{b}(x-\tau_{1},y-\tau_{2})\overline{W_{A_{2}}^{\bm{\mu}}}(y,\tau_{2})d\tau_{1}d\tau_{2}\bigg]\bm{\nu}.
\nonumber
\end{eqnarray}
Taking the QLCT of the first bracket of Eq. $(\ref{huc1})$ with Theorem  \ref{le43} and the linear property of the QLCT,  we have
$$e^{-\bm{\mu}(\frac{d_{1}}{2b_{1}}u^2+\frac{d_{2}}{2b_{2}}v^2)}
\bigg[\mathcal{L}^{\bm{\mu},\bm{\mu}}[f_{a}](u,v)\mathcal{L}^{\bm{\mu},\bm{\mu}}[g_{a}](u,v)
-\mathcal{L}^{\bm{\mu},\bm{\mu}}[f_{b}](u,v)\mathcal{L}^{\bm{\mu},\bm{\mu}}[\overline{ g_{b}}](u,v)\bigg].
$$
By direct computations, the QLCT  of the second  bracket of the equation $(\ref{huc1})$ equals to
$$e^{-\bm{\mu}(\frac{d_{1}}{2b_{1}}u^2-\frac{d_{2}}{2b_{2}}v^2)}\bigg[\mathcal{L}^{\bm{\mu},-\bm{\mu}}{[f_{b}]}(u,v)
\mathcal{L}^{\bm{\mu},-\bm{\mu}}[\overline{g}_{a}](u,v)
+\mathcal{L}^{\bm{\mu},-\bm{\mu}}[f_{a}](u,v)\mathcal{L}^{\bm{\mu},-\bm{\mu}}{[g_{b}]}(u,v)\bigg]\bm{\nu}$$
 then it completes the proof of equation $(\ref{hu83})$.  The proof of the Eq.$(\ref{h83})$  follows in an analogous way, we omit it.
\end{nproof}
\begin{remark}\label{re47}
When  $ A_{1}=\left(
             \begin{array}{cc}
               \cos \alpha & \sin \alpha \\
               -\sin \alpha  & \cos\alpha  \\
             \end{array}
           \right)
, A_{2}=\left(
             \begin{array}{cc}
               \cos \beta & \sin \beta \\
               -\sin \beta & \cos\beta  \\
             \end{array}
           \right)$,
we obtain the spatial convolution theorems of the quaternion fractional Fourier transform (QFrFT). Compared to the convolution theorems of the QFrFT  in \cite{guanlei2008fractional}, our spatial  convolution operation is more concise and  easier  to be implemented.  From Eqs. (\ref{h81})  and   (\ref{h83}), the  spatial convolution theorem maintains the exact product in the QLCT domain without the chirp multipliers. On the other hand, the functions in \cite{guanlei2008fractional}  are multiplied three times by the different chirp signals.  It is hard to realize because in communication systems it is nearly impossible to generate a chirp signal accurately \cite{wei2012new}. In our spatial  convolution structure, two chirp signals are used. It should be pointed out that, however,  this spatial convolution operation can not keep the simplicity in the $\mathbb{H}$ space, as the QLCT  with parameter matrices
$A_{i}=\left(
             \begin{array}{cc}
               a_{i} &b_{i} \\
               c_{i}  &d_{i}  \\
             \end{array}
           \right), i=1,2  $ of  the  convolution
of two quaternionic functions  is equal to  the summation of products of the QLCT  with parameter matrices
$\tilde{A}_{i} =\left(
             \begin{array}{cc}
               a_{i} &b_{i} \\
              \frac{ c_{i}}{2} -\frac{1}{2b_{i}} &\frac{ d_{i}}{2}  \\
             \end{array}
           \right)$, $ i=1,2 $ of  their components.
\end{remark}

\begin{corollary}(\textbf{Product theorem on $ \mathbb{H}$})\label{th48}
Suppose that $f,g\in L^{1}\bigcap L^{2}(\mathbb{R}^2,\mathbb{H}),$  and $\mathcal{L}^{\bm{\mu},\bm{\mu}}[f]\in L^{p}(\mathbb{R}^2,\mathbb{H}),\mathcal{L}^{\bm{\mu},\bm{\mu}}[g]\in L^{1}(\mathbb{R}^2,\mathbb{H}),$
then
\begin{eqnarray}\label{hu90}
&&\mathcal{L}^{\bm{\mu},\bm{\mu}}[ e^{\bm{\mu}\frac{a_{1}}{2b_{1}}x^{2}} fg(x,y)e^{\bm{\mu}\frac{a_{2}}{2b_{2}}y^{2}}](u,v)\\
=&&\left  (  \mathcal{L}^{\bm{\mu},\bm{\mu}} [ f_{a}]\circledast_{A_{1}^{-1},A_{2}^{-1}} \mathcal{L}^{\bm{\mu},\bm{\mu}}[ g_{a}](u,v)-\mathcal{L}^{\bm{\mu},\bm{\mu}}  [ f_{b}]\circledast_{A_{1}^{-1},A_{2}^{-1}}\mathcal{L} [ \overline{ g_{b}}](u,v)\right)\nonumber\\
&&+\left (  \mathcal{L}^{\bm{\mu},-\bm{\mu}} [  f_{b}]\circledast_{A_{1}^{-1},A_{2}^{-1}}\mathcal{L}^{\bm{\mu},-\bm{\mu}} [ \overline{ g_{a}} ](u,v)+\mathcal{L}^{\bm{\mu},-\bm{\mu}}  (  f_{a})\circledast_{A_{1}^{-1},A_{2}^{-1}}\mathcal{L}^{\bm{\mu},-\bm{\mu}} [ g_{b}](u,v)\right)\bm{\nu}, \nonumber
\end{eqnarray}
and 
\begin{eqnarray}\label{h90}
&&\mathcal{L}_{ \breve{A}_{1},  \breve{A}_{2}}^{\bm{\mu},\bm{\mu}}  [ fg](u,v)\\
&=&\left  (  \mathcal{L}^{\bm{\mu},\bm{\mu}} [ f_{a}]\circledast_{A_{1}^{-1},A_{2}^{-1}}\mathcal{L}^{\bm{\mu},\bm{\mu}}[  g_{a}](u,v)-\mathcal{L}^{\bm{\mu},\bm{\mu}}  [ f_{b}]\circledast_{A_{1}^{-1},A_{2}^{-1}}\mathcal{L}^{\bm{\mu},\bm{\mu}} [ \overline{ g_{b}}](u,v)\right)\nonumber\\
&&+\left ( \mathcal{L}^{\bm{\mu},-\bm{\mu}} [  f_{b}]\circledast_{A_{1}^{-1},A_{2}^{-1}}\mathcal{L}^{\bm{\mu},-\bm{\mu}} [ \overline{ g_{a}} ](u,v)+\mathcal{L}^{\bm{\mu},-\bm{\mu}}  [f_{a}]\circledast_{A_{1}^{-1},A_{2}^{-1}}\mathcal{L}^{\bm{\mu},-\bm{\mu}} [  g_{b}](u,v)\right)\bm{\nu}.\nonumber
\end{eqnarray}

\end{corollary}
\begin{nproof}
Since
\begin{eqnarray*}
&&\mathcal{L}^{\bm{\mu},\bm{\mu}} [ e^{\bm{\mu}\frac{a_{1}}{2b_{1}}x^{2}}fg(x,y)e^{\bm{\mu}\frac{a_{2}}{2b_{2}}y^{2}}]\\
=&&\mathcal{L}^{\bm{\mu},\bm{\mu}}[  e^{\bm{\mu}\frac{a_{1}}{2b_{1}}x^{2}}(f_{a}g_{a}-f_{b}\overline{ g_{b}})e^{\bm{\mu}\frac{a_{2}}{2b_{2}}y^{2}}]
+\mathcal{L}^{\bm{\mu},\bm{\mu}} [  e^{\bm{\mu}\frac{a_{1}}{2b_{1}}x^{2}}(f_{b}\overline{g_{a}}\bm{\nu}+f_{a}g_{b}\bm{\nu})e^{\bm{\mu}\frac{a_{2}}{2b_{2}}y^{2}}],
\end{eqnarray*}
using Theorem \ref{le44}, we   obtain
\begin{eqnarray*}
\mathcal{L} ^{\bm{\mu},\bm{\mu}}[  e^{\bm{\mu}\frac{a_{1}}{2b_{1}}x^{2}}(f_{a}g_{a}-f_{b}\overline{ g_{b}})e^{\bm{\mu}\frac{a_{2}}{2b_{2}}y^{2}}](u,v)
=  \mathcal{L}^{\bm{\mu},\bm{\mu}} [  f_{a}]\circledast_{A_{1}^{-1},A_{2}^{-1}}\mathcal{L}^{\bm{\mu},\bm{\mu}} [  g_{a}]-\mathcal{L}^{\bm{\mu},\bm{\mu}}[ f_{b}]\circledast_{A_{1}^{-1},A_{2}^{-1}}\mathcal{L}^{\bm{\mu},\bm{\mu}} [ \overline{ g_{b}}](u,v).
\end{eqnarray*}
Then we only need to prove the following equation.
 \begin{eqnarray*}
&&\mathcal{L}^{\bm{\mu},\bm{\mu}} (  e^{\bm{\mu}\frac{a_{1}}{2b_{1}}x^{2}}e^{-\bm{\mu}\frac{a_{2}}{2b_{2}}y^{2}}[f_{b}\overline{g_{a}}\bm{\nu}+f_{b}g_{b}\bm{\nu}])(u,v)\\
=&&\left  ( \mathcal{L}^{\bm{\mu},-\bm{\mu}}   [ f_{b}]\circledast_{A_{1}^{-1},A_{2}^{-1}}\mathcal{L}^{\bm{\mu},-\bm{\mu}}   [  \overline{g_{a}}](u,v)+\mathcal{L}^{\bm{\mu},-\bm{\mu}}  [ f_{b}] \circledast_{A_{1}^{-1},A_{2}^{-1}}\mathcal{L}^{\bm{\mu},-\bm{\mu}}  [ g_{b} ](u,v)\right)\bm{\nu}.
\end{eqnarray*}
Since
 \begin{eqnarray*}
&&\mathcal{L}_{A_{1}^{-1},A_{2}^{-1}}^{\bm{\mu},\bm{\mu}}[  \mathcal{L}^{\bm{\mu},-\bm{\mu}}  [  f_{b}]\circledast_{A_{1}^{-1},A_{2}^{-1}} \mathcal{L}^{\bm{\mu},-\bm{\mu}}  [  \overline{g_{a}}](u,v)\bm{\nu}](x,y)\\
&&=\int_{-\infty}^{\infty}\int_{-\infty}^{\infty}\frac{1}{\sqrt{-2\bm{\mu} b_{1}\pi}}e^{\bm{\mu}(-\frac{d_{1}}{2b_{1}}u^{2}+\frac{1}{b_{1}}ux-\frac{a_{1}}{2b_{1}}x^{2})}\\
&& \int_{-\infty}^{\infty}\int_{-\infty}^{\infty}\frac{1}{\sqrt{-2\bm{\mu} b_{1}\pi}}e^{\bm{\mu}\tau_{1}(\tau_{1}-u)(-\frac{d_{1}}{b_{1}})}\mathcal{L}^{\bm{\mu},-\bm{\mu}} [f_{b}]( u-\tau_{1} ,v-\tau_{2})
\mathcal{L}^{\bm{\mu},-\bm{\mu}} [\overline{g_{a}}](\tau_{1} ,\tau_{2})\\
&&\frac{1}{\sqrt{2\bm{\mu} b_{2}\pi}}e^{\bm{\mu}\tau_{2}(\tau_{2}-v)(\frac{d_{2}}{b_{2}})}\bm{\nu} d\tau_{1}d\tau_{2}
 \frac{1}{\sqrt{-2\bm{\mu} b_{2}\pi}}e^{\bm{\mu}(-\frac{d_{2}}{2b_{2}}v^{2}+\frac{1}{b_{2}}vy-\frac{a_{2}}{2b_{2}}y^{2})}dudv,
 \end{eqnarray*}
replacing  $ u-\tau_{1} ,v-\tau_{2}$ by $s,t,$ respectively. We have
\begin{eqnarray*}
&&\int_{-\infty}^{\infty}\int_{-\infty}^{\infty}\frac{1}{\sqrt{-2\bm{\mu} b_{1}\pi}}e^{\bm{\mu}(-\frac{d_{1}}{2b_{1}}(s+\tau_{1})^{2}+\frac{1}{b_{1}}( s+\tau_{1})x-\frac{a_{1}}{2b_{1}}x^{2})}\\
&& \int_{-\infty}^{\infty}\int_{-\infty}^{\infty}\frac{1}{\sqrt{-2\bm{\mu} b_{1}\pi}}e^{\bm{\mu}\tau_{1}s(\frac{d_{1}}{b_{1}})}\mathcal{L}^{\bm{\mu},-\bm{\mu}} [f_{b}](s ,t)\mathcal{L}^{\bm{\mu},-\bm{\mu}} [\overline{g_{a}}   ](\tau_{1} ,\tau_{2})\frac{1}{\sqrt{2\bm{\mu} b_{2}\pi}}e^{-\bm{\mu}\tau_{2}t(\frac{d_{2}}{b_{2}})}\bm{\nu} d\tau_{1}d\tau_{2}\\
&&  \frac{1}{\sqrt{-2\bm{\mu} b_{2}\pi}}e^{\bm{\mu}(-\frac{d_{2}}{2b_{2}}( t+\tau_{2} )^{2}+\frac{1}{b_{2}}( t+\tau_{2} )y-\frac{a_{2}}{2b_{2}}y^{2})}dsdt\\
=&&\int_{-\infty}^{\infty}\int_{-\infty}^{\infty}\frac{1}{\sqrt{-2\bm{\mu} b_{1}\pi}}e^{\bm{\mu}(-\frac{d_{1}}{2b_{1}}(s+\tau_{1})^{2}+\frac{1}{b_{1}}( s+\tau_{1})x-\frac{a_{1}}{2b_{1}}x^{2})}\\
&&  \int_{-\infty}^{\infty}\int_{-\infty}^{\infty}\frac{1}{\sqrt{-2\bm{\mu} b_{1}\pi}}e^{\bm{\mu}\tau_{1}s(\frac{d_{1}}{b_{1}})}\mathcal{L}^{\bm{\mu},-\bm{\mu}} [f_{b}](s ,t)\mathcal{L}^{\bm{\mu},-\bm{\mu}} [\overline{g_{a}}](\tau_{1} ,\tau_{2})\frac{1}{\sqrt{2\bm{\mu} b_{2}\pi}} e^{-\bm{\mu}\tau_{2}t(\frac{d_{2}}{b_{2}})}d\tau_{1}d\tau_{2}\\
&&  \frac{1}{\sqrt{2\bm{\mu} b_{2}\pi}}e^{\bm{\mu}(\frac{d_{2}}{2b_{2}}( t+\tau_{2} )^{2}-\frac{1}{b_{2}}( t+\tau_{2} )y+\frac{a_{2}}{2b_{2}}y^{2})}\bm{\nu} dsdt\\
=&&e^{\bm{\mu}\frac{a_{1}}{2b_{1}}x^{2}}e^{-\bm{\mu}\frac{a_{2}}{2b_{2}}y^{2}}f_{b}\overline{g_{a}}(x,y)\bm{\nu}.
\end{eqnarray*}
With the analogous computations, we can obtain the following equation,
 \begin{eqnarray*}
\mathcal{L}_{\A_{1}^{-1},\A_{2}^{-1}}^{\bm{\mu},\bm{\mu}} [\mathcal{L}^{\bm{\mu},-\bm{\mu}} [  f_{a}]\circledast_{A_{1}^{-1},A_{2}^{-1}}\mathcal{L}^{\bm{\mu},-\bm{\mu}}  [  g_{b}]\bm{\nu}](x,y)
=e^{\bm{\mu}\frac{a_{1}}{2b_{1}}x^{2}}e^{-\bm{\mu}\frac{a_{2}}{2b_{2}}y^{2}}f_{a}g_{b}(x,y)\bm{\nu}.
\end{eqnarray*}
The proof of Eq.  (\ref{hu90})  is now completed. By using the analogous argument of Eq. (\ref{hu90}), Eq.(\ref{h90}) can be derived accordingly.
\end{nproof}

\subsection{Applications}
 In this subsection, we primarily investigate the applications of the  QLCT and its spatial  convolution operator   $\circledast_{A_{1},A_{2}}$,
 including  defining the correlation operation for the  QLCT,  solving integral equations, and partial differential equations, designing multiplicative filters.

\subsubsection{The correlation theorem}
The classical correlation operator $ \circledcirc$ is defined by the  2-D  classical  convolution operator $\ast$ as follows.
Given two integrable complex-valued  functions,
$$ f\circledcirc g(x,y) := \overline{f(-\cdot,-\cdot)}\ast g(x,y),$$
taking the 2-D Fourier transform of both sides gives the cross-correlation theorem,
$$ \widehat{(f\circledcirc g)}(u,v)=\overline{\hat{f}(u,v)} \hat{g}(u,v).$$
When $f=g$,  then the correlation theorem reduces to the Wiener-Khinchin theorem, which is the foundation theorem for
random signal processing. The quaternionic correlation is  defined
in  S. J. Sangwine's paper  \cite{sangwine1999hypercomplex}.  The Wiener-Khinchin theorem for discrete quaternion Fourier transform is proven in \cite{ell2000hypercomplex}, which not only allows the computation of vector correlations for color images but also gives the theoretical foundation for 3-D and 4-D    quaternionic random signal processing\cite{took2011augmented}, such as color image processing \cite{ell2000hypercomplex},   simulation and estimation problems \cite{navarro2016semi,ginzberg2013quaternion}, wind modelling \cite{gou2015three}.  The QLCT is more flexible than the QFT and  QFRFT, therefore we plan to extend the Wiener-Khinchin theorem to the 2-D QLCT domain. But, first, we need to define the correlation operator for the QLCT.
With  the spatial  convolution operator $\circledast_{A_{1},A_{2}} $, we can define the correlation operator  $\circledcirc_{A_{1},A_{2}}$ for the QLCT as follows.

\begin{definition}\label{def61}
Suppose  $f\in L^{p}(\mathbb{R}^{2}, \mathbb{H}), g \in L^{1}(\mathbb{R}^{2}, \mathbb{H}),$ then the correlation operator $\circledcirc_{A_{1},A_{2}} $ is defined by
\begin{eqnarray}\label{hu61}
f \circledcirc_{A_{1},A_{2}}  g(x,y)&:=&\overline{f(-\cdot,-\cdot)}\circledast_{A_{1},A_{2}}  g(x,y)\nonumber\\
&=&\int_{\mathbb{R }^{2}}W_{A_{1}}^{\bm{\mu}}(x,\tau_{1} )\overline{f(-\tau _{1},-\tau _{2})}
g(x-\tau _{1},y-\tau _{2})W_{A_{2}}^{\bm{\mu}}(y,\tau _{2})d\tau _{1}d\tau _{2}.\nonumber\\
\end{eqnarray}
\end{definition}
\begin{remark}
As a result, this is a reasonable approach to defining the correlation operation for the QLCT, since
when $a_{i}=0,i=1,2$, the Formula   (\ref{hu61})  can be reduced to the classical case of the  2-D Fourier transform and QFT in  Todd A.Ell etc's book \cite{ell2014quaternion}.

\end{remark}


\begin{theorem}(\textbf{Correlation theorem} on $  \mathbb{H}(\bm{\mu}$))\label{th62}
Suppose $f\in L^{p}(\mathbb{R}^{2}, \mathbb{H}(\bm{\mu}))$ and $g\in L^{1}(\mathbb{R}^{2}, \mathbb{H}(\bm{\mu})),$
then
\begin{eqnarray*}
\mathcal{L}^{\bm{\mu},\bm{\mu}}{[(f\circledcirc_{A_{1},A_{2}} g)]}(u,v)&&=e^{-\bm{\mu}\frac{d_{1}}{2b_{1}}u^{2}}e^{-\bm{\mu}\frac{d_{2}}{2b_{2}}v^{2}}
\mathcal{L}^{\bm{\mu},\bm{\mu}}{[\bar{f}]}(-u,-v)\mathcal{L}^{\bm{\mu},\bm{\mu}}{[g]}(u,v).\\
\mathcal{L}^{\bm{\mu},\bm{\mu}}{[(f\circledcirc_{A_{1},A_{2} }g)]}(u,v)&&=
\mathcal{L}_{\tilde{A}_{1},\tilde{A}_{2}}^{\bm{\mu},\bm{\mu}}[\bar{f}](-u,-v)\mathcal{L}_{\tilde{A}_{1},\tilde{A}_{2}}^{\bm{\mu},\bm{\mu}}[g](u,v ).
\end{eqnarray*}
\end{theorem}
\begin{nproof}
Using the analogous argument as in the proof of Theorem \ref{le43}, we can easily carry out the proof of this theorem.
\end{nproof}
\begin{remark}
when $a_{i}= d_{i}=0,b_{i}=1,i=1,2$, then $\mathcal{L}_{\tilde{A}_{1},\tilde{A}_{2}}^{\bm{\mu},\bm{\mu}}[\bar{f}](-u,-v)=
\overline{\mathcal{L}_{\tilde{A}_{1},\tilde{A}_{2}}^{\bm{\mu},\bm{\mu}}[f](u,v)}$,
Correlation Theorem $\ref{th62}$ reduces to the Wiener-Khinchin theorem  of the 2-D Fourier transform.
\end{remark}

\begin{theorem}(\textbf{Correlation theorem} on $ \mathbb{H}$)\label{th63}
For two given quaternionic functions $f\in L^{p}(\mathbb{R}^{2},\mathbb{H}), g\in L^{1}(\mathbb{R}^{2},\mathbb{H}),$
$f=f_{a}+\bm{\nu}\overline{ f_{b}},$
$g=g_{a}+\bm{\nu}\overline{ g_{b}},$
$ f_{a},f_{b},g_{a},g_{b} \in \mathbb{H}(\bm{\mu})$ which is defined in Eq.$(\ref{ab1})$.
then
\begin{eqnarray*}
&&\mathcal{L}^{\bm{\mu},\bm{\mu}}[(f\circledcirc_{A_{1},A_{2}} g)](u,v)\nonumber \\
=&&e^{-\bm{\mu}(\frac{d_{1}}{2b_{1}}u^2+\frac{d_{2}}{2b_{2}}v^2)}
\big(\mathcal{L}^{\bm{\mu},\bm{\mu}}[\bar{f}_{a}](-u,-v)\mathcal{L}^{\bm{\mu},\bm{\mu}}[g
_{a}](u,v)
+\mathcal{L}^{\bm{\mu},\bm{\mu}}(f_{b})(-u,-v)\mathcal{L}^{\bm{\mu},\bm{\mu}}[\overline{g}_{b}](u,v)\big)
\\&&-e^{-\bm{\mu}(\frac{d_{1}}{2b_{1}}u^2-\frac{d_{2}}{2b_{2}}v^2)}
\big(\mathcal{L}^{\bm{\mu},-\bm{\mu}}[f_{b}](-u,-v)\mathcal{L}^{\bm{\mu},-\bm{\mu}}[\overline{g}_{a}](u,v)
\\&&-\mathcal{L}^{\bm{\mu},-\bm{\mu}}[\bar{f}_{a}](-u,-v)\mathcal{L}^{\bm{\mu},-\bm{\mu}}[b_{b}](u,v)\big)\bm{\nu}.
\end{eqnarray*}
\begin{eqnarray*}
&&\mathcal{L}^{\bm{\mu},\bm{\mu}}[(f\circledcirc_{A_{1},A_{2}} g)](u,v)\nonumber \\
=&&
\big(\mathcal{L}_{\tilde{A}_{1},\tilde{A}_{2}}^{\bm{\mu},\bm{\mu}}[\bar{f}_{a}](-u,-v)
\mathcal{L}_{\tilde{A}_{1},\tilde{A}_{2}}^{\bm{\mu},\bm{\mu}}[g
_{a}](u,v)
+\mathcal{L}_{\tilde{A}_{1},\tilde{A}_{2}}^{\bm{\mu},\bm{\mu}}[f_{b}](-u,-v)
\mathcal{L}_{\tilde{A}_{1},\tilde{A}_{2}}^{\bm{\mu},\bm{\mu}}[\overline{g}_{b}](u,v)\big)\\
&&-\big(\mathcal{L}^{\bm{\mu},-\bm{\mu}}_{\tilde{A}_{1},\tilde{A}_{2}}[f_{b}](-u,-v)
\mathcal{L}^{\bm{\mu},-\bm{\mu}}_{\tilde{A}_{1},\tilde{A}_{2}}[\overline{g}_{a}](u,v)
-\mathcal{L}^{\bm{\mu},-\bm{\mu}}_{\tilde{A}_{1},\tilde{A}_{2}}[\bar{f}_{a}](-u,-v)
\mathcal{L}^{\bm{\mu},-\bm{\mu}}_{\tilde{A}_{1},\tilde{A}_{2}}[b_{b}](u,v)\big)\bm{\nu}.
\end{eqnarray*}

\end{theorem}
\begin{nproof}
It can be proved using the same method as Theorem \ref{th46}.
\end{nproof}

 \subsubsection{Solving integral equations}
Many problems in engineering and mechanics can be converted  into 2-D Fredholm integral equations,
such as 2-D  heat conduction equations with the Cauchy problem.
In terms of using the potential theorem, a 3-D Laplace equation with boundary conditions can be transformed into a  2-D boundary integral equation with a  weakly singular kernel function \cite{kress1989linear}. The deblurring of 2-D images can be modelled as a 2-D integral equation with a smooth kernel function \cite{chan1996conjugate,rajan2003simultaneous}.

We are interested in the kernels involving trigonometric functions in Fredholm integral equations \cite{polyanin2008handbook}. By using the Convolution theorem \ref{th46}, we can solve the following Fredholm integral equation of the first kind. Let

\begin{eqnarray}\label{h1051}
 \int_{\mathbb{R}^{2} } K(x, y, \tau_{1}, \tau_{2})f(\tau_{1}, \tau_{2} )d\tau_{1}d\tau_{2}=g(x, y),
\end{eqnarray}
where $  K(x, y, \tau_{1}, \tau_{2})=W_{A_{1}}^{\bm{\mu}}(x,\tau_{1} )r(x-\tau _{1},y-\tau _{2})W_{A_{2}}^{\bm{\mu}}(y,\tau _{2})$. If $r, g \in L^{1}(\mathbb{R}^{2}, \mathbb{H}),$ then a solution $f \in L^{1}(\mathbb{R}^{2}, \mathbb{H}(\bm{\mu})) $ is desired.

According to Definition \ref{def41}, the left hand side of the Eq.(\ref{h1051}) is $r \circledast_{A_{1},A_{2}} f(x,y)  \in L^{1}(\mathbb{R}^{2}, \mathbb{H}).$ Applying the QLCT to both sides of  the Eq.(\ref{h1051}), we have
\begin{eqnarray*}
&&e^{-\bm{\mu}(\frac{d_{1}}{2b_{1}}u^2+\frac{d_{2}}{2b_{2}}v^2)}
[\mathcal{L}^{\bm{\mu},\bm{\mu}}[r_{a}](u,v)\mathcal{L}^{\bm{\mu},\bm{\mu}}[f](u,v)]\\
&&+e^{-\bm{\mu}(\frac{d_{1}}{2b_{1}}u^2-\frac{d_{2}}{2b_{2}}v^2)}
[\mathcal{L}^{\bm{\mu},-\bm{\mu}}[r_{b}](u,v)\mathcal{L}^{\bm{\mu},-\bm{\mu}}[\overline{f}](u,v)]\bm{\nu}\\
=&&\mathcal{L}^{\bm{\mu},\bm{\mu}}[g_{a}](u,v) +\mathcal{L}^{\bm{\mu},-\bm{\mu}}[g_{b}](u,v) \bm{\nu}.
\end{eqnarray*}
If $ \mathcal{L}^{\bm{\mu},\bm{\mu}}[r_{a}](u,v) \neq 0$, then we have
\begin{eqnarray}\label{h1053}
\mathcal{L}^{\bm{\mu},\bm{\mu}}[f](u,v)
=\left(e^{\bm{-\mu}\left(\frac{d_{1}}{2b_{1}}u^{2}+\frac{d_{2}}{2b_{2}}v^{2}\right)}\mathcal{L}^{\bm{\mu},\bm{\mu}}[r_{a}](u,v)\right)^{-1}\mathcal{L}^{\bm{\mu},\bm{\mu}}[g_{a}](u,v).
\end{eqnarray}
If the right hand side of Eq.   (\ref{h1053})  is a function of $ L^{1}(\mathbb{R}^{2}, \mathbb{H}(\bm{\mu})),$ then using the inverse QLCT, we can solve the Eq.  (\ref{h1053}) as
\begin{eqnarray*}
f(x,y)
=\int_{\mathbb{R}^{2}}\frac{1}{\sqrt{-2\pi b_{1}\bm{\mu}}} e^{\bm{\mu}(\frac{1}{b_{1}}xu-\frac{a_{1}}{2b_{1}}x^{2})} \left((\mathcal{L}^{\bm{\mu},\bm{\mu}}[r_{a}](u,v))^{-1}\mathcal{L}^{\bm{\mu},\bm{\mu}}[g_{a}](u,v)
\right)\frac{1}{\sqrt{-2\pi b_{2}\bm{\mu}}} e^{\bm{\mu}(\frac{1}{b_{2}}yv-\frac{a_{2}}{2b_{2}}y^{2})}dudv.
\end{eqnarray*}
If $ \mathcal{L}^{\bm{\mu},-\bm{\mu}}[r_{b}](u,v) \neq 0$, then
 \begin{eqnarray*}
\bar{f}(x,y)
=\int_{\mathbb{R}^{2}}\frac{1}{\sqrt{-2\pi b_{1}\bm{\mu}}} e^{\bm{\mu}(\frac{1}{b_{1}}xu-\frac{a_{1}}{2b_{1}}x^{2})} \left((\mathcal{L}^{\bm{\mu},-\bm{\mu}}[r_{b}](u,v))^{-1}\mathcal{L}^{\bm{\mu},-\bm{\mu}}[g_{b}](u,v)
\right)\frac{1}{\sqrt{2\pi b_{2}\bm{\mu}}} e^{-\bm{\mu}(\frac{1}{b_{2}}yv-\frac{a_{2}}{2b_{2}}y^{2})}dudv.
\end{eqnarray*}

\subsubsection{Solving partial differential  equations (PDEs)}
In this subsection,  let $\mathbb{S}(\mathbb{R}^{2}, \mathbb{ H})$ denote the Schwartz space from $ \mathbb{R}^2$ into $ \mathbb{H},$ or space of rapidly decreasing smooth $( C^{\infty})$ quaternionic functions on $ \mathbb{R}^{2}$. That is, 
 \begin{eqnarray*}
\mathbb{S}(\mathbb{R}^{2}, \mathbb{ H}):=\left\{ f\in C^{\infty} : \sup_{(x,y)\in \mathbb{R}^2}(1+(x^{2}+y^{2}))^{k/2}| \partial^{\alpha}f(x, y)|< \infty \right\},
 \end{eqnarray*}
 where $k$ is a non-negative integer, $ \alpha =(\alpha_{1}, \alpha_{2} )$ is a multi-index of non-negatives $\alpha_{1}$ and $\alpha_{2},$ $ \partial^{ \alpha} f:=\frac{\partial^{\alpha_{1}+\alpha_{2}}}{\partial x^{\alpha_{1}}\partial y^{\alpha_{2}}} f$, $ C^{\infty}$  is the set of smooth functions from  $ \mathbb{R}^2$ into $ \mathbb{H}.$
 To proceed, we need the following  useful lemma which gives the second order partial and mixed derivative of $f$.
The proof can be followed by the integration by parts.

\begin{lemma} \label{dlct2}
Suppose $f\in \mathbb{S}(\mathbb{R}^{2}, \mathbb{ H})$, then
$$\mathcal{L}^{\bm{\mu},\bm{\mu}}\left[\frac{\partial^{2}f(x,y) }{\partial x \partial y}\right](u,v)
=\int_{\mathbb{R}^{2}}\bm{\mu} (\frac{a_{1}}{b_{1}}x-\frac{1}{b_{1}}u)K_{A_{1}}^{\bm{\mu}}(x,u)f(x,y)
K_{A_{2}}^{\bm{\mu}}(y,v)(\frac{a_{2}}{b_{2}}y-\frac{1}{b_{2}}v)\bm{\mu}dxdy$$

and
\begin{eqnarray*}
\mathcal{L}^{\bm{\mu},\bm{\mu}}\left[\frac{\partial^{2}f(x,y) }{\partial x^{2}  }\right](u,v)
=&&\int_{\mathbb{R}^{2}}\bm{\mu}\frac{a_{1}}{b_{1}}K_{A_{1}}^{\bm{\mu}}(x,u)f(x,y)
K_{A_{2}}^{\bm{\mu}}(y,v)dxdy\\
&&-\int_{\mathbb{R}^{2}}\frac{a_{1}^{2}}{b_{1}^{2}}x^{2}K_{A_{1}}^{\bm{\mu}}(x,u)f(x,y)
K_{A_{2}}^{\bm{\mu}}(y,v)dxdy\\
&&+\int_{\mathbb{R}^{2}}2\frac{a_{1}}{b_{1}^{2}}xuK_{A_{1}}^{\bm{\mu}}(x,u)f(x,y)
K_{A_{2}}^{\bm{\mu}}(y,v)dxdy\\
&&-\int_{\mathbb{R}^{2}}\frac{1}{b_{1}^{2}}u^{2}K_{A_{1}}^{\bm{\mu}}(x,u)f(x,y)
K_{A_{2}}^{\bm{\mu}}(y,v)dxdy.
\end{eqnarray*}

\end{lemma}

\begin{example}
Solve the system of second order  real linear  partial differential equations,
\begin{eqnarray} \label{h1016}
   \left\{
 \begin{array}{llll}
    \frac{\partial^{2} f_{0}(x,y)}{\partial x \partial y} -\frac{a_{2}}{b_{2}}y\frac{\partial f_{1}(x,y)}{\partial x}
-\frac{a_{1}}{b_{1}}x\frac{\partial f_{1}(x,y)}{\partial y} -\frac{a_{1}}{b_{1}}\frac{a_{2}}{b_{2}}xyf_{0}(x,y)=g_{0}(x,y),
 &  \\[1.5ex]
    \frac{\partial^{2} f_{1}(x,y)}{\partial x \partial y} +\frac{a_{2}}{b_{2}}y\frac{\partial f_{0}(x,y)}{\partial x}
+\frac{a_{1}}{b_{1}}x\frac{\partial f_{0}(x,y)}{\partial y} -\frac{a_{1}}{b_{1}}\frac{a_{2}}{b_{2}}xyf_{1}(x,y)=g_{1}(x,y),&
 &  \\[1.5ex]
    \frac{\partial^{2} f_{2}(x,y)}{\partial x \partial y} +\frac{a_{2}}{b_{2}}y\frac{\partial f_{3}(x,y)}{\partial x}
-\frac{a_{1}}{b_{1}}x\frac{\partial f_{3}(x,y)}{\partial y} +\frac{a_{1}}{b_{1}}\frac{a_{2}}{b_{2}}xyf_{2}(x,y)=g_{2}(x,y),&
 &  \\[1.5ex]
    \frac{\partial^{2} f_{3}(x,y)}{\partial x \partial y} -\frac{a_{2}}{b_{2}}y\frac{\partial f_{2}(x,y)}{\partial x}
+\frac{a_{1}}{b_{1}}x\frac{\partial f_{1}(x,y)}{\partial y} +\frac{a_{1}}{b_{1}}\frac{a_{2}}{b_{2}}xyf_{3}(x,y)=g_{3}(x,y).&
 &
 \end{array}\right.
\end{eqnarray}

The system of equations (\ref{h1016}) is first converted  into the following second-order quaternion  partial differential equation, 
\begin{eqnarray}\label{h108}
\frac{\partial^2 f(x,y)}{\partial x \partial y} +\frac{a_{2}}{b_{2}}y\frac{\partial f(x,y)}{\partial x}\bm{\mu}
+\bm{\mu}\frac{a_{1}}{b_{1}}x\frac{\partial f(x,y)}{\partial y} +\bm{\mu}\frac{a_{1}}{b_{1}}\frac{a_{2}}{b_{2}}xyf(x,y)\bm{\mu}=g(x,y)
\end{eqnarray}
where
\begin{eqnarray*}
 f(x,y)=f_{0}(x,y)+ \bm{\mu}f_{1}(x,y)+\bm{\nu}f_{2}(x,y)+\bm{\eta}f_{3}(x,y),\\
 g(x,y)=g_{0}(x,y)+ \bm{\mu}g_{1}(x,y)+\bm{\nu}g_{2}(x,y)+\bm{\eta}g_{3}(x,y).
\end{eqnarray*}
If $f(x,y)$  and $g(x,y)$ belong to $\mathbb{S}(\mathbb{R}^{2}, \mathbb{ H}),$ then taking  the QLCT 
of both sides, we have
\begin{eqnarray*}
\bm{\mu}\frac{1}{b_{1}}\frac{1}{b_{2}}uv \mathcal{L}^{\bm{\mu},\bm{\mu}}[f](u,v)\bm{\mu}=\mathcal{L}^{\bm{\mu},\bm{\mu}}[g](u,v).
\end{eqnarray*}
Finally, if $ \frac{1}{uv}\mathcal{L}^{\bm{\mu},\bm{\mu}}[g](u,v)\in L^{1}(\mathbb{R}^{2},\mathbb{H}),$ then 
we obtain the solution of the Eq.$(\ref{h108})$ as following,
\begin{eqnarray}\label{h109}
f(x,y)&=&\int_{\mathbb{R}^{2}}\frac{1}{\sqrt{-2\pi \bm{\mu}b_{1}}}e^{
 \bm{\mu}(-\frac{d_{1}}{2b_{1}}u^{2}+\frac{1}{b_{1}}xu-\frac{a_{1}}{2b_{1}}x^{2})} \nonumber\\ && \bm{\mu}\frac{b_{1}b_{2}}{uv}\mathcal{L}^{\bm{\mu},\bm{\mu}}[g](u,v)\bm{\mu}\frac{1}{\sqrt{-2\pi\bm{\mu}b_{2}}}e^{\bm{\mu}(-\frac{d_{2}}{2b_{2}}v^{2}+\frac{1}{b_{2}}yv-\frac{a_{2}}{2b_{2}}y^{2})}dudv.
\end{eqnarray}

Suppose  $\mathcal{L}^{\bm{\mu},\bm{\mu}}[g](u,v)=uv e^{-(u^2+v^2)},$ then
\begin{eqnarray*}
\mathcal{L}^{\bm{\mu},\bm{\mu}}[f](u,v)=-b_{1}b_{2}e^{-(u^2+v^2)},
\end{eqnarray*}
according to   Eq. (\ref{h109}), we  obtain the solution
\begin{eqnarray*}
f(x,y)&=&-b_{1}b_{2}\sqrt{\frac{2b_{1}\pi}{2b_{1}-d_{1}\bm{\mu}}}\sqrt{\frac{2b_{2}\pi}{2b_{2}-d_{2}\bm{\mu}}}\\
&& e^{((4b_{1}-2d_{1}\bm{\mu})a_{1}+2/(2b_{1}(4b_{1}+2d_{1}\bm{\mu}))x^{2}}
e^{((4b_{2}-2d_{2}\bm{\mu})a_{2}+2/(2b_{2}(4b_{2}+2d_{2}\bm{\mu}))y^{2}}.
\end{eqnarray*}

\end{example}

\begin{example}
Solve the linear second order  elliptic partial differential equation
\begin{eqnarray}\label{h1081}
\frac{\partial^2 f(x,y)}{\partial x^2} +\frac{\partial^2 f(x,y)}{\partial y^2}+
2\bm{\mu}\frac{a_{2}}{b_{2}}y\frac{\partial f(x,y)}{\partial y}
+2\bm{\mu}\frac{a_{1}}{b_{1}}x\frac{\partial f(x,y)}{\partial x} \nonumber \\
+\left(2\bm{\mu}\frac{a_{1}}{b_{1}}+2\bm{\mu}\frac{a_{2}}{b_{2}}-\frac{a_{1}^2}{b_{1}^2}x^2-\frac{a_{2}^2}{b_{2}^2}y^2\right)f(x,y)=g(x,y).
\end{eqnarray}
 If $f, g$ belong to  $\mathbb{S}(\mathbb{R}^2, \mathbb{H}(\bm{\mu}))$, taking the  QLCT
of both sides, then we can get an algebraic equation:
\begin{eqnarray*}
\left(\bm{\mu}\frac{a_{1}}{b_{1}}+\bm{\mu}\frac{a_{2}}{b_{2}}\right)\mathcal{L}^{\bm{\mu},\bm{\mu}}[f](u,v)
-\left(\frac{1}{b_{1}^2}u^2+\frac{1}{b_{2}^2}v^2\right)\mathcal{L}^{\bm{\mu},\bm{\mu}}[f](u,v)
=\mathcal{L}^{\bm{\mu},\bm{\mu}}[g](u,v).
\end{eqnarray*}
$ \Longleftrightarrow$
\begin{eqnarray*}
\mathcal{L}^{\bm{\mu},\bm{\mu}}[f](u,v)
=\left(\bm{\mu}\frac{a_{1}}{b_{1}}+\bm{\mu}\frac{a_{2}}{b_{2}}-\frac{1}{b_{1}^2}u^2-\frac{1}{b_{2}^2}v^2\right)^{-1}\mathcal{L}^{\bm{\mu},\bm{\mu}}[g](u,v).
\end{eqnarray*}
Set $\mathcal{L}^{\bm{\mu},\bm{\mu}}[r](u,v)
=\left(\bm{\mu}\frac{a_{1}}{b_{1}}+\bm{\mu}\frac{a_{2}}{b_{2}}-\frac{1}{b_{1}^2}u^2-\frac{1}{b_{2}^2}v^2\right)^{-1}e^{\bm{\mu}(\frac{d_{1}}{2b_{1}}u^2+\frac{d_{2}}{2b_{2}}v^2)} $,  it is easy to verify that  $\mathcal{L}^{\bm{\mu},\bm{\mu}}[r]\in L^{2}(\mathbb{R}^{2}, \mathbb{H}(\bm{\mu}))$.
 According to Lemma \ref{L1} and Theorem \ref{le43}, we can have the solution of PDE  (\ref{h1081})  as
$$ f(x,y)=\int_{\mathbb{R}^{2}}W_{A_{1}}^{\bm{\mu}}(x,\tau_{1} )r(\tau _{1},\tau _{2})g(x-\tau _{1},y-\tau _{2})W_{A_{2}}^{\bm{\mu}}(y,\tau_{2} )d\tau _{1}d\tau _{2}. $$

\end{example}
\begin{remark}
When $ \bm{\mu}=\i$, which can be considered as the complex imaginary unit, then $ \mathbb{H}(\i)$ can be regarded as the complex field $ \mathbb{C}$. If   $f(x,y)$ and $ g(x,y)$  both belong to $ \mathbb{H}(\i)$, then
the partial differential equation $(\ref{h1081})$ reduces to the complex  elliptic PDE \cite{florian2001functional},
\begin{eqnarray*}
\frac{\partial^2 f(x,y)}{\partial x^2} +\frac{\partial^2 f(x,y)}{\partial y^2}+
2\i\frac{a_{2}}{b_{2}}y\frac{\partial f(x,y)}{\partial y}
+2\i\frac{a_{1}}{b_{1}}x\frac{\partial f(x,y)}{\partial x}\\
+\left(2\i\frac{a_{1}}{b_{1}}+2\i\frac{a_{2}}{b_{2}}-\frac{a_{1}^2}{b_{1}^2}x^2-\frac{a_{2}^2}{b_{2}^2}y^2\right)f(x,y)=g(x,y).
\end{eqnarray*}
\end{remark}
\begin{example}
Solve the special case of the system of  second order real  partial differential equations \cite{smirnov1964course} about  anisotropic elastic media.
\begin{eqnarray}\label{h1015}
   \left\{
 \begin{array}{llll}
   &&\frac{\partial^2 f_{r}(x,y)}{\partial x^2} +\frac{\partial^2 f_{r}(x,y)}{\partial y^2}+
\frac{\partial^2 f_{r}(x,y)}{\partial x \partial y}
-\left(\frac{a_{2}}{b_{2}}y+2\frac{a_{1}}{b_{1}}x\right)\frac{\partial f_{i}(x,y)}{\partial x}
-\left(\frac{a_{1}}{b_{1}}x+2\frac{a_{2}}{b_{2}}y\right)\frac{\partial f_{i}(x,y)}{\partial y} \\
&&+\left(-2\frac{a_{1}}{b_{1}}-2\frac{a_{2}}{b_{2}}\right)f_{i}(x,y)+\left(-\frac{a_{1}^2}{b_{1}^2}x^2-\frac{a_{2}^2}{b_{2}^2}y^2-\frac{a_{1}}{b_{1}}\frac{a_{2}}{b_{2}}xy\right)f_{r}(x,y)
=g_{r}(x,y), &  \\[1.5ex]
 &&\frac{\partial^2 f_{i}(x,y)}{\partial x^2} +\frac{\partial^2 f_{i}(x,y)}{\partial y^2}+
\frac{\partial^2 f_{i}(x,y)}{\partial x \partial y}
+\left(\frac{a_{2}}{b_{2}}y+2\frac{a_{1}}{b_{1}}x\right)\frac{\partial f_{r}(x,y)}{\partial x}
+\left(\frac{a_{1}}{b_{1}}x+2\frac{a_{2}}{b_{2}}y\right)\frac{\partial f_{r}(x,y)}{\partial y} \\
&&+\left(2\frac{a_{1}}{b_{1}}+2\frac{a_{2}}{b_{2}})f_{r}(x,y)+(-\frac{a_{1}^2}{b_{1}^2}x^2-\frac{a_{2}^2}{b_{2}^2}y^2-\frac{a_{1}}{b_{1}}\frac{a_{2}}{b_{2}}xy\right)f_{i}(x,y)
=g_{i}(x,y).&
 \end{array}\right.
\end{eqnarray}

 Firstly, utilizing $ f(x,y)=f_{r}(x,y)+ \bm{\mu}f_{i}(x,y),$
 $g(x,y)=g_{r}(x,y)+ \bm{\mu}g_{i}(x,y),$
 we can  convert the  PDEs $(\ref{h1015})$ into the following  second-order quaternion partial differential equation
\begin{eqnarray}\label{h10152}
&&\frac{\partial^2 f(x,y)}{\partial x^2} +\frac{\partial^2 f(x,y)}{\partial y^2}+
\frac{\partial^2 f(x,y)}{\partial x \partial y} \nonumber\\
&&+\bm{\mu}\left(\frac{a_{2}}{b_{2}}y+2\frac{a_{1}}{b_{1}}x\right)\frac{\partial f(x,y)}{\partial x}
+\bm{\mu}\left(\frac{a_{1}}{b_{1}}x+2\frac{a_{2}}{b_{2}}y\right)\frac{\partial f(x,y)}{\partial y} \nonumber\\
&&+\left(2\bm{\mu}\frac{a_{1}}{b_{1}}+2\bm{\mu}\frac{a_{2}}{b_{2}}-\frac{a_{1}^2}{b_{1}^2}x^2-\frac{a_{2}^2}{b_{2}^2}y^2-\frac{a_{1}}{b_{1}}\frac{a_{2}}{b_{2}}xy\right)f(x,y)
=g(x,y).
\end{eqnarray}
If $f, g$ belong to  $\mathbb{S}(\mathbb{R}^2, \mathbb{H}(\bm{\mu})),$  applying  the QLCT
to Eq. (\ref{h10152}), then we can obtain an algebraic equation below:
\begin{eqnarray*}
\left(\bm{\mu}\frac{a_{1}}{b_{1}}+\bm{\mu}\frac{a_{2}}{b_{2}}\right)\mathcal{L}^{\bm{\mu},\bm{\mu}}[f](u,v)
-\left(\frac{1}{b_{1}^2}u^2+\frac{1}{b_{1}}\frac{1}{b_{2}}uv+\frac{1}{b_{2}^2}v^2\right)\mathcal{L}^{\bm{\mu},\bm{\mu}}[f](u,v)
=\mathcal{L}^{\bm{\mu},\bm{\mu}}[g](u,v).
\end{eqnarray*}
$\Longleftrightarrow$
\begin{eqnarray*}
\mathcal{L}^{\bm{\mu},\bm{\mu}}[f](u,v)
=
\left(\bm{\mu}\frac{a_{1}}{b_{1}}+\bm{\mu}\frac{a_{2}}{b_{2}}-\frac{1}{b_{1}^2}u^2-\frac{1}{b_{1}}\frac{1}{b_{2}}uv-\frac{1}{b_{2}^2}v^2\right)^{-1}\mathcal{L}^{\bm{\mu},\bm{\mu}}[g](u,v).
\end{eqnarray*}
Set $\mathcal{L}^{\bm{\mu},\bm{\mu}}[r](u,v)=
\left(\bm{\mu}\frac{a_{1}}{b_{1}}+\bm{\mu}\frac{a_{2}}{b_{2}}-\frac{1}{b_{1}^2}u^2-\frac{1}{b_{1}}\frac{1}{b_{2}}uv-\frac{1}{b_{2}^2}v^2\right)^{-1}e^{\bm{\mu}(\frac{d_{1}}{2b_{1}}u^2+\frac{d_{2}}{2b_{2}}v^2)}$
,   it is easy to verify that  $\mathcal{L}^{\bm{\mu},\bm{\mu}}[r] \in L^{2}(\mathbb{R}^{2}, \mathbb{H(\bm{\mu})})$.
 According to Lemma \ref{L1} and Theorem \ref{le43}, we can have the solution of PDEs  $(\ref{h1015})$ as following,
$$ f(x,y)=\int_{\mathbb{R}^{2}}W_{\A_{1}}^{\bm{\mu}}(x,\tau_{1} )r(\tau _{1},\tau _{2})g(x-\tau _{1},y-\tau _{2})W_{\A_{2}}^{\bm{\mu}}(y,\tau_{2} )d\tau _{1}d\tau _{2}. $$

\end{example}

 \subsubsection{Multiplicative filters}

In this subsection,  the applications of proposed theorems in designing  multiplicative filters  are investigated.

\begin{example}
The multiplicative filters in the QLCT domain are shown in Figs \ref{ff1} and \ref{ff2} for quaternionic functions taking value in $\mathbb{H}(\bm{\mu})$.  If the parameter matrices of  QLCT of input and output quaternionic functions are the same as
$A_{i}=\left(
\begin{array}{cc}
	a_{i} &b_{i} \\
	c_{i}  &d_{i}  \\
\end{array}
\right), i=1,2.$
According to Eq.(\ref{hu81}) of the Theorem \ref{le43},   the transfer function is
$e^{-\bm{\mu}(\frac{d_{1}}{2b_{1}}u^{2}+\frac{d_{2}}{2b_{2}}v^{2})}\mathcal{L}^{\bm{\mu},\bm{\mu}}{[g]}(u,v)$  in  the multiplicative filter as in Fig.\ref{ff1}.
Otherwise, according to Eq.(\ref{h81}) of the Theorem \ref{le43},  the transfer function is  $ \mathcal{L}^{\bm{\mu},\bm{\mu}}[g](u,v)$ as in Fig. \ref{ff2}. Under this case,  the output function  is  $f_{out}(x,y)$ which  has the QLCT  with
parameter matrices
$\hat{A}_{i}=\left(
             \begin{array}{cc}
               a_{i} &b_{i} \\
               2c_{i}+1/b_{i}  &2d_{i}  \\
             \end{array}
           \right), i=1,2 $.
The different types of multiplicative filters, such as low pass, high pass, band pass, band stop, and so on, can be obtained by designing the transfer functions.
\end{example}

\begin{figure}
  \begin{center}\label{ff1}
    \includegraphics[width=10cm]{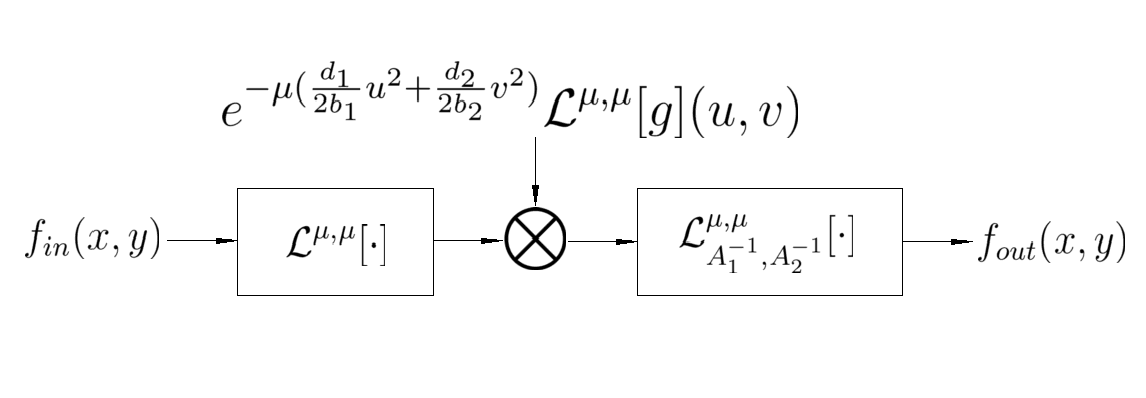}
     \caption{  Multiplicative filter  for  QLCT domain with same  parameter matrixes.  }
     \label{ff1}
  \end{center}
\end{figure}

\begin{figure}
  \begin{center}\label{ff2}
    \includegraphics[width=10cm]{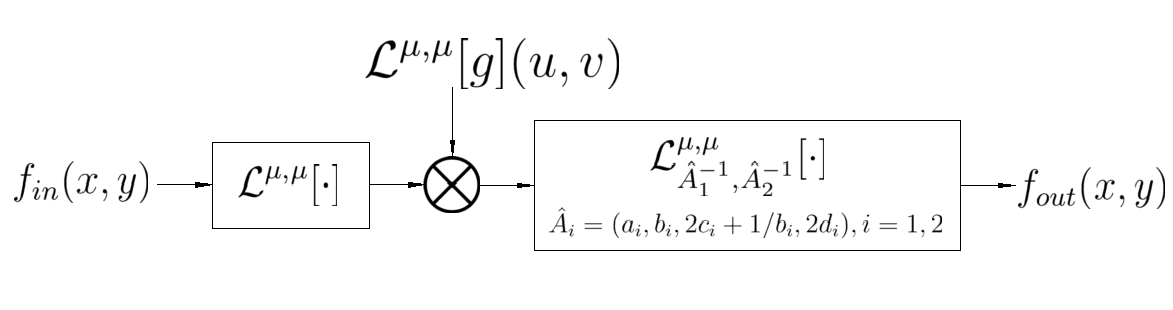}
     \caption{  Multiplicative filter for  QLCT domain with different  parameter matrixes.  }
     \label{ff2}
 \end{center}
\end{figure}

\section{Spectral convolution theorem}\label{covsec4}

The purpose of this section is to introduce the second type of convolution operator, motivated by the spectral domain representation, namely the spectral convolution operator of the QLCT. It preserves the convolution theorem for the classical Fourier transform. That is, the  QLCT of a convolution of two quaternionic functions is the pointwise product of their QLCTs. 

\subsection{Spectral convolution $ \bigstar $}
The spectral convolution operator is defined in the QLCT domain. Precisely,
\begin{definition} (\textbf{Spectral convolution operator $ \bigstar$})
	Let $f\in L^{2}(\mathbb{R}^{2}, \mathbb{H}),$ and $g\in L^{2}(\mathbb{R}^{2}, \mathbb{H}),$  the spectral convolution  operator $ \bigstar$ of the  QLCT is defined by
	\begin{eqnarray} \label{Type2}
		[f\bigstar g ](u,v) :=\mathcal{L}^{\bm{\mu},\bm{\mu}}_{A_1^{-1},A_2^{-1}} [\mathcal{L}^{\bm{\mu},\bm{\mu}}[f]\mathcal{L}^{\bm{\mu},\bm{\mu}}[g]](u,v).
	\end{eqnarray} 
	It can also be represented as 	
	\begin{eqnarray}\label{c2}
	f\bigstar g(x,y)&=&
		e^{\bm{\mu}(\frac{d_{1}}{2b_{1}}u^{2}+\frac{d_{2}}{2b_{2}}v^{2})}
		\int_{\mathbb{R}^{2} }W_{A_{1}}^{\bm{\mu}}(x,\tau_{1} )f_{a}(x-\tau_{1},y-\tau_{2})g_{a}(\tau_{1},\tau_{2})W_{A_{2}}^{\bm{\mu}}(y,\tau_{2} )d\tau_{1}d\tau_{2}\nonumber\\
		&&+e^{\bm{\mu}(-\frac{d_{1}}{2b_{1}}u^{2}-\frac{d_{2}}{2b_{2}}v^{2})}
		\int_{\mathbb{R}^{2} }E_{A_{1}}^{\bm{\mu}}(x,\tau_{1} )f_{b}(x-\tau_{1},y-\tau_{2})\bm{\nu}g_{a}(-\tau_{1},\tau_{2})W_{A_{2}}^{\bm{\mu}}(y,\tau_{2} )d\tau_{1}d\tau_{2}\nonumber\\
		&&+e^{\bm{\mu}(\frac{d_{1}}{2b_{1}}u^{2}+\frac{d_{2}}{2b_{2}}v^{2})}
		\int_{\mathbb{R}^{2} }W_{A_{1}}^{\bm{\mu}}(x,\tau_{1} )f_{a}(x-\tau_{1},-y+\tau_{2})\bm{\nu}\bar{g}_{b}(\tau_{1},\tau_{2})W_{\check{A}_{2}}^{\bm{\mu}}(\tau_{2},y)d\tau_{1}d\tau_{2}\nonumber\\
		&&-e^{\bm{\mu}(-\frac{d_{1}}{2b_{1}}u^{2}-\frac{d_{2}}{2b_{2}}v^{2})}
		\int_{\mathbb{R}^{2} }E_{A_{1}}^{\bm{\mu}}(x,\tau_{1} )f_{b}(x-\tau_{1},-y+\tau_{2})\bar{g}_{b}(-\tau_{1},\tau_{2})\overline{W_{A_{2}}^{\bm{\mu}}}(\tau_{2},y )d\tau_{1}d\tau_{2},
	\end{eqnarray}
	where
	\begin{eqnarray*}
		E_{A_{1}}^{\bm{\mu}}(x,\tau_{1}):=\frac{1}{\sqrt{-2\pi b_{1}\bm{\mu}}}e^{\bm{\mu}(-\tau_{1}x\frac{a_{1}}{b_{1}})},\quad
		\check{A}_{2}:=\left(
		\begin{array}{cc}
			-a_{2} &b_{2} \\
			c_{2} &-d_{2}  \\
		\end{array}
		\right).
	\end{eqnarray*}
\end{definition}
\begin{remark}
	Suppose  $f\in L^{p}(\mathbb{R}^{2}, \mathbb{H}),$ and $g\in L^{1}(\mathbb{R}^{2}, \mathbb{H}),$  then the expression of $ \bigstar$ in (\ref{c2}) is also well defined.
\end{remark}

\begin{theorem} (\textbf{Spectral convolution theorem in $\mathbb{H}$})\label{h0312}
If $f\in L^{p}(\mathbb{R}^{2}, \mathbb{H}),$ and $g\in L^{1}(\mathbb{R}^{2}, \mathbb{H}),$ then
\begin{eqnarray*}
\mathcal{L}^{\bm{\mu},\bm{\mu}}[f\bigstar g](u,v) =\mathcal{L}^{\bm{\mu},\bm{\mu}}[f](u,v)\mathcal{L}^{\bm{\mu},\bm{\mu}}[g](u,v).
\end{eqnarray*}
\end{theorem}

\begin{nproof}
Since  $f$ and $g\in L^{1}(\mathbb{R}^{2}, \mathbb{H}),$ then  their components $ f_{a}, f_{b}, g_{a}, g_{b}$ all belong to
$L^{1}(\mathbb{R}^{2}, \mathbb{H}).$ According to  Lemma \ref{le42}, it follows that  $\mathcal{L}^{\bm{\mu},\bm{\mu}}[f\bigstar g(x,y)](u,v)$ is well defined.
\par
Let $ I_{n}, n=1,2,3,4$ denote the four components of Eq. (\ref{c2}), respectively. By direct calculation, we have
 \begin{eqnarray*}
 \mathcal{L}^{\bm{\mu},\bm{\mu}}[f\bigstar ](u,v)
 &=&\mathcal{L}^{\bm{\mu},\bm{\mu}}[I_{1}+I_{2}+I_{3}+I_{4}](u,v)\\
 &=& \mathcal{L}^{\bm{\mu},\bm{\mu}}[f_{a}](u,v)\mathcal{L}^{\bm{\mu},\bm{\mu}}[g_{a}](u,v)
 +\mathcal{L}^{\bm{\mu},-\bm{\mu}}[f_{b}](u,v)\bm{\nu}\mathcal{L}^{\bm{\mu},\bm{\mu}}[g_{a}](u,v)\\
 &&+\mathcal{L}^{\bm{\mu},\bm{\mu}}[f_{a}](u,v)\mathcal{L}^{\bm{\mu},-\bm{\mu}}[g_{b}](u,v)\bm{\nu}
 -\mathcal{L}^{\bm{\mu},-\bm{\mu}}[f_{b}](u,v)\mathcal{L}^{-\bm{\mu},\bm{\mu}}[\bar{g}_{b}](u,v)\\
  &=& \mathcal{L}^{\bm{\mu},\bm{\mu}}[f_{a}+f_{b}\bm{\nu}](u,v) \mathcal{L}^{\bm{\mu},\bm{\mu}}[g_{a}+g_{b}\bm{\nu}](u,v),
 \end{eqnarray*}
which completes the proof.
\end{nproof}

\begin{remark}
    The spectral convolution  structure of quaternionic functions is a natural generalization of the convolution structure associated with the 2-D  Fourier transform, hence it would be convenient to achieve the multiplicative filter in the QLCT domain. When $a_{i}=d_{i}=0, b_{i}=1, i=1,2$,Theorem 8 is equivalent to the convolution theorem of the two-sided QFT.
\end{remark}
\begin{corollary}
If $f\in L^{p}(\mathbb{R}^{2}, \mathbb{H}),$ and $g\in L^{1}(\mathbb{R}^{2}, \mathbb{H}),$ then
\begin{eqnarray} \label{h1011}
\mathcal{F}[f\bigstar g](u,v)=\mathcal{F}[f](u,v)\mathcal{F}[g](u,v),
\end{eqnarray}
\noindent where $ \mathcal{F}[\cdot]$ denotes the two-sided QFT.
\end{corollary}
\begin{remark}
Compared to the convolution theorem in \cite{bahri2013convolution},  the spectral convolution operator $\bigstar $ holds the excellent property of classical convolution associated with the  2-D Fourier transform. Meanwhile, the proposed general convolution operation is more succinct compared to the convolution theorem in \cite{pei2001efficient}, which is called spectrum-product quaternion convolution.
\end{remark}

If the parameter matrices  $ A_1$ and $A_2 $ of the  QLCT are replaced by  $ \A_{\alpha }:=\left(
             \begin{array}{cc}
               \cos \alpha & \sin \alpha \\
               -\sin \alpha  & \cos\alpha  \\
             \end{array}
           \right)$ and 
$ \A_{\beta}=\left(
             \begin{array}{cc}
               \cos \beta & \sin \beta \\
               -\sin \beta & \cos\beta  \\
             \end{array}
           \right)$,  respectively, then the spectral convolution Theorem $\ref{h0312}$ of the QLCT
           reduces to the convolution formula in quaternion fractional Fourier transform as follows.

\begin{corollary}
	If $f\in L^{p}(\mathbb{R}^{2}, \mathbb{H}),$ and $g\in L^{1}(\mathbb{R}^{2}, \mathbb{H}),$ then
\begin{eqnarray*}
\mathcal{L}_{\A_{\alpha}, \A_{\beta} }^{\bm{\mu},\bm{\mu}}[f\bigstar g ](u,v)=\mathcal{L}_{\A_{\alpha}, \A_{\beta}  }^{\bm{\mu},\bm{\mu}}[f](u,v)\mathcal{L}_{\A_{\alpha}, \A_{\beta} }^{\bm{\mu},\bm{\mu}}[g](u,v).
\end{eqnarray*}

\end{corollary}
\subsection{Applications}

\subsubsection{The correlation theorem}
\begin{definition}  \label{def261}
For  $f\in L^{p}(\mathbb{R}^{2}, \mathbb{H}), g \in L^{1}(\mathbb{R}^{2}, \mathbb{H}),$   the correlation operator $ \circledR$  is defined  by
\begin{eqnarray}\label{hu261}
f  \circledR g(x,y) &:=& (f(\cdot,\cdot) \bigstar \overline{g(-\cdot,-\cdot)})(x,y)\nonumber\\
&=&e^{\bm{\mu}(\frac{d_{1}}{2b_{1}}u^{2}+\frac{d_{2}}{2b_{2}}v^{2})}
\int_{\mathbb{R}^{2} }W_{A_{1}}^{\bm{\mu}}(x,\tau_{1} )f_{a}(x-\tau_{1},y-\tau_{2})\overline{g_{a}}(-\tau_{1},-\tau_{2})W_{A_{2}}^{\bm{\mu}}(y,\tau_{2} )d\tau_{1}d\tau_{2}\nonumber\\
&&+e^{\bm{\mu}(-\frac{d_{1}}{2b_{1}}u^{2}-\frac{d_{2}}{2b_{2}}v^{2})}
\int_{\mathbb{R}^{2} }E_{A_{1}}^{\bm{\mu}}(x,\tau_{1} )f_{b}(x-\tau_{1},y-\tau_{2})\bm{\nu}\overline{g_{a}}(\tau_{1},-\tau_{2})W_{A_{2}}^{\bm{\mu}}(y,\tau_{2} )d\tau_{1}d\tau_{2}\nonumber\\
&&-e^{\bm{\mu}(\frac{d_{1}}{2b_{1}}u^{2}+\frac{d_{2}}{2b_{2}}v^{2})}
\int_{\mathbb{R}^{2} }W_{A_{1}}^{\bm{\mu}}(x,\tau_{1} )f_{a}(x-\tau_{1},-y+\tau_{2})\bm{\nu}\bar{g}_{b}(-\tau_{1},-\tau_{2})W_{\check{A}_{2}}^{\bm{\mu}}(\tau_{2},y)d\tau_{1}d\tau_{2}\nonumber\\
&&+e^{\bm{\mu}(-\frac{d_{1}}{2b_{1}}u^{2}-\frac{d_{2}}{2b_{2}}v^{2})}
\int_{\mathbb{R}^{2} }E_{A_{1}}^{\bm{\mu}}(x,\tau_{1} )f_{b}(x-\tau_{1},-y+\tau_{2})\bar{g}_{b}(\tau_{1},-\tau_{2})\overline{W_{A_{2}}^{\bm{\mu}}}(\tau_{2},y )d\tau_{1}d\tau_{2}.\nonumber\\
\end{eqnarray}
\end{definition}

Applying the QLCT on both sides of above, we obtain the correlation theorem.
\begin{theorem} $($\textbf{Correlation theorem  of $ \circledR$}$)$  \label{th262}
If $f$ and $g\in L^{1}(\mathbb{R}^{2}, \mathbb{H}),$ then 
\begin{eqnarray*}
\mathcal{L}^{\bm{\mu},\bm{\mu}}[f  \circledR g](u,v)=\mathcal{L}^{\bm{\mu},\bm{\mu}}[f](u,v)\mathcal{L}^{\bm{\mu},\bm{\mu}}[\overline{g}](-u,-v).
\end{eqnarray*}
\end{theorem}

\begin{nproof}
Let $ I_{n}, n=1,2,3,4$ denote the four components of Eq.$(\ref{hu261})$, respectively. By direct calculations, we have
 \begin{eqnarray*}
 \mathcal{L}^{\bm{\mu},\bm{\mu}}[f\circledR g](u,v)
 &=&\mathcal{L}^{\bm{\mu},\bm{\mu}}[I_{1}+I_{2}+I_{3}+I_{4}](u,v)\\
 &=& \mathcal{L}^{\bm{\mu},\bm{\mu}}[f_{a}](u,v)\mathcal{L}^{\bm{\mu},\bm{\mu}}[\overline{g_{a}}](-u,-v)
 +\mathcal{L}^{\bm{\mu},-\bm{\mu}}[f_{b}](u,v)\bm{\nu}\mathcal{L}^{\bm{\mu},\bm{\mu}}[\overline{g_{a}}](u,v)\\
 &&-\mathcal{L}^{\bm{\mu},\bm{\mu}}[f_{a}](u,v)\mathcal{L}^{\bm{\mu},-\bm{\mu}}[g_{b}](u,v)\bm{\nu}
 +\mathcal{L}^{\bm{\mu},-\bm{\mu}}[f_{b}](u,v)\mathcal{L}^{-\bm{\mu},\bm{\mu}}[\bar{g}_{b}](u,v)\\
  &=& \mathcal{L}^{\bm{\mu},\bm{\mu}}[f_{a}+f_{b}\bm{\nu}](u,v) \mathcal{L}^{\bm{\mu},\bm{\mu}}[\overline{g_{a}}-g_{b}\bm{\nu}](-u,-v)\\
  &=& \mathcal{L}^{\bm{\mu},\bm{\mu}}[f](u,v) \mathcal{L}^{\bm{\mu},\bm{\mu}}[\overline{g}](-u,-v).
 \end{eqnarray*}
\end{nproof}

\begin{remark}
\begin{enumerate}
  \item When $a_{i}= d_{i}=0,b_{i}=1,i=1,2$, we have $\mathcal{L}^{\bm{\mu},\bm{\mu}}[\bar{g}](-u,-v)=\overline{\mathcal{L}^{\bm{\mu},\bm{\mu}}[g](u,v)}$,
that is, for $f\in L^{1}(\mathbb{R}^{2}, \mathbb{H}),$ then
\begin{eqnarray*}
 \mathcal{L}^{\bm{\mu},\bm{\mu}}[f\circledR f](u,v)=|\mathcal{L}^{\bm{\mu},\bm{\mu}}[f](u,v) |^{2}.
 \end{eqnarray*}
The correlation theorem is completely parallel to the 2-D Fourier transform, as is evident.
  \item
As shown in Theorems $\ref{th62}$, $\ref{th63}$, and $\ref{th262}$  convolution operations  $\circledast_{A_{1},A_{2}} $ and $\bigstar$ can implement correlation operations of two quaternionic functions.

\end{enumerate}
\end{remark}

 \subsubsection{Solving partial differential equations}
\begin{example}
Solve the  second order  quaternion  partial differential equation
\begin{eqnarray}\label{hu14}
&&\frac{\partial^2 f(x,y)}{\partial x^2}+ \frac{\partial^2 f(x,y)}{\partial y^2}+\frac{a_{1}^{2}}{b_{1}^{2}}x^{2}f(x,y)
-f(x,y)\bm{\mu}\frac{a_{2}}{b_{2}} +\frac{a_{2}^{2}}{b_{2}^{2}}y^{2}f(x,y)+\nonumber \\
&&\sqrt{\bm{\mu}2\pi b_1 }\bm{\mu} \frac{2a_{1}}{b_{1}} p_{1}\circledast_{A_1,A_2} xf(x,y)\sqrt{\bm{\mu}2\pi b_2 }+ \sqrt{\bm{\mu}2\pi b_1 }\frac{2a_{2}}{b_{2}} yf\circledast_{A_1,A_2} p_{2}(x,y)\bm{\mu}\sqrt{\bm{\mu}2\pi b_2 }\nonumber\\
=&&g(x,y), 
\end{eqnarray}
where $p_{1}(x,y)=e^{-\bm{\mu}\frac{a_{1}x^{2}}{2b_{1}}}\delta'(x)\delta(y),
p_{2}(x,y)=e^{-\bm{\mu}\frac{a_{2}y^{2}}{2b_{2}}}\delta'(y)\delta(x), a_1 \neq 0;
 f,g \in \mathbb{S}(\mathbb{R}^2, \mathbb{H}).$
\end{example}
 Due to $\mathcal{L}^{\bm{\mu},\bm{\mu}}[p_{1}](u,v)=\bm{\mu}\frac{1}{b_{1}} u e^{\bm{\mu}(\frac{d_{1}u^{2}}{2b_{1}}+\frac{d_{2}v^{2}}{2b_{2}})},\mathcal{L}^{\bm{\mu},\bm{\mu}}[p_{2}](u,v)=\bm{\mu}\frac{1}{b_{2}} v e^{\bm{\mu}(\frac{d_{1}u^{2}}{2b_{1}}+\frac{d_{2}v^{2}}{2b_{2}})},$
 applying the QLCT on both sides of  Eq. (\ref{hu14}),  and by Lemma  \ref{dlct2}  and Convolution Theorem $\ref{le43},$ we have
 \begin{eqnarray*}
\bm{\mu}\frac{a_{1}}{b_{1}}\mathcal{L}^{\bm{\mu},\bm{\mu}}[f](u,v)-\frac{1}{b_{1}^{2}}u^{2}\mathcal{L}^{\bm{\mu},\bm{\mu}}[f](u,v)
-\frac{1}{b_{2}^{2}}v^{2}\mathcal{L}^{\bm{\mu},\bm{\mu}}[f](u,v)
=\mathcal{L}^{\bm{\mu},\bm{\mu}}[g](u,v).
\end{eqnarray*}
It follows that
\begin{eqnarray*}
\bigg(\bm{\mu}\frac{a_{1}}{b_{1}}-\frac{1}{b_{1}^{2}}u^{2}-\frac{1}{b_{2}^{2}}v^{2}\bigg)\mathcal{L}^{\bm{\mu},\bm{\mu}}[f](u,v)
&=&\mathcal{L}^{\bm{\mu},\bm{\mu}}[g](u,v)\\
\mathcal{L}^{\bm{\mu},\bm{\mu}}[f](u,v)&=&    \bigg(\bm{\mu}\frac{a_{1}}{b_{1}}-\frac{1}{b_{1}^{2}}u^{2}-\frac{1}{b_{2}^{2}}v^{2}\bigg)^{-1} \mathcal{L}^{\bm{\mu},\bm{\mu}}[g](u,v).
\end{eqnarray*}
since $ a_1 \neq 0$, then $\bigg(\bm{\mu}\frac{a_{1}}{b_{1}}-\frac{1}{b_{1}^{2}}(u)^{2}-\frac{1}{b_{2}^{2}}(v)^{2}\bigg)^{-1} $ is well defined.
By the spectral convolution theorem \ref{h0312},  and let $ l(x,y) := \mathcal{L}_{A_{1}^{-1},A_{2}^{-1}}^{\bm{\mu},\bm{\mu}}\left[ \bigg(\bm{\mu}\frac{a_{1}}{b_{1}}-\frac{1}{b_{1}^{2}}(\cdot)^{2}-\frac{1}{b_{2}^{2}}(\cdot)^{2}\bigg)^{-1}\right](x,y),$ we obtain the solution of second order quaternion partial differential equation
\begin{eqnarray*}
f(x,y)&=& l\bigstar g(x,y).
\end{eqnarray*}
\subsubsection{Multiplicative filter}
\begin{figure}\label{Fig222}
  \begin{center}
    \includegraphics[width=10cm]{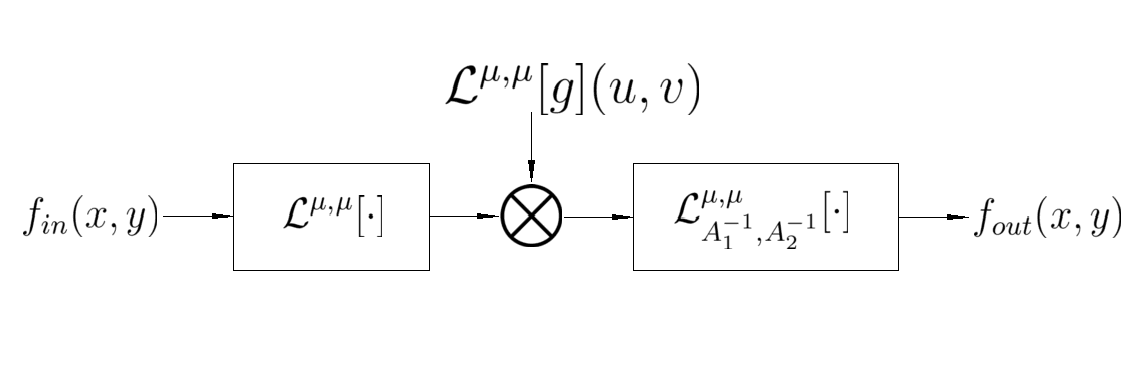}
     \caption{  Multiplicative filter for  QLCT domain  .  }
 \label{f3}
   \end{center}
\end{figure}

According to Theorem  $\ref{h0312}$ of quaternionic functions,
the model of the multiplicative filter in the QLCT domain is shown in Fig.\ref{f3}.
Let $ \mathcal{L}^{\bm{\mu},\bm{\mu}}[g](u,v)$ be the transformed function.
 We want to design a lowpass  filter by QLCT with the passband of $ \{ (u,v)| u_{1}<u< u_{2}, v_{1}<v< v_{2},\}$. We can design the lowpass filter as in the conventional case
          $$g(x,y)=\mathcal{L}_{A_{1}^{-1},A_{2}^{-1}}^{\bm{\mu},\bm{\mu}}[\mathcal{L}^{\bm{\mu},\bm{\mu}}[g]](x,y),$$
where
\begin{eqnarray*}
 \mathcal{L}^{\bm{\mu},\bm{\mu}}[g](u,v)
  =\left\{
 \begin{array}{llll}
    1, &  \quad  u_{1}<u< u_{2}, v_{1}<v< v_{2},\\[1.5ex]
    0, &  \quad \text{otherwise}.
     \end{array}\right.
\end{eqnarray*}
Accordingly, the QLCT of the quaternionic  output $f_{out}(x,y)$ and the quaternionic input $f_{in}(x,y)$ have the following relationship
\begin{eqnarray*}
    \mathcal{L}^{\bm{\mu},\bm{\mu}}[f_{out}](u,v)=\left\{
 \begin{array}{llll}
      \mathcal{L}^{\bm{\mu},\bm{\mu}}[f_{in}](u,v), &  \quad  u_{1}<u< u_{2}, v_{1}<v< v_{2},\\[1.5ex]
0, & \text{otherwise}.
 \end{array}\right.
\end{eqnarray*}

\begin{example}
 Suppose the original size of the given figure is $ N\times M $. This example is about two lowpass filters, namely $ H_1$ and $H_2$, they are given as follows
\begin{eqnarray}\label{H1}
  H_1(u,v)=\left\{
 \begin{array}{llll}
    1, &  \quad  N/4<u< 3N/4, M/4<v< 3M/4,\\[1.5ex]
    0, &\quad \text{otherwise},
     \end{array}\right.
\end{eqnarray}

\begin{eqnarray*}
  H_2(u,v)=\left\{
 \begin{array}{llll}
    1, &  \quad  N/6<u< 3N/6, M/6<v< 3M/6,\\[1.5ex]
    0, &\quad \text{otherwise}.
     \end{array}\right.
\end{eqnarray*}
Lowpass filters $ H_1$ and  $ H_2$ have domain areas of $\frac{1}{ 4}NM$ and $\frac{1}{ 9}NM$, respectively. As a result, lowpass filter $ H_1$ receives more information than lowpass filter $ H_2$.
Two well-known objective image quality metrics, namely the peak-signal-to-noise ratio (PSNR) and signal-to-noise (SNR) are analyzed. The  PSNR is a ratio between the maximum power of the signal and the power of corrupting noise.  In the case of a reference image $f$ and a test image $g$ of size $M\times N$, the PSNR between $f$ and $g$ is defined as follows:
      $$ \mbox{PSNR}(f,g) :=10 \log_{10}\left(\frac{M \times N}{ \mbox{MSN}(f,g)}\right), $$
      where $ \mbox{MSN}(f,g) :=\frac{1}{M \times N} \sum^{M}_{i=1}\sum^{N}_{j=1}(f_{ij}-g_{ij})^{2}.$
    When the value of MSN approaches zero, the PSNR value approaches infinity,  it shows that the higher value of PSNR implies  the higher image quality.
    \begin{table}[h!]
\caption{PSNR and SNR comparison values for the house image }\label{tab2}
\centering
\begin{tabular}{|c|c|c|}
  \hline
   $A_1$=(0,2,-1/2,2);
   $A_2$=(0,2,-1/2,4) & SNR& PSNR  \\
       \hline
Image (B) & 7.6181     &  3.0019 \\
    \hline
 Image (C) &  4.8974 & 0.2812    \\
    \hline

\end{tabular}
\end{table}

\begin{figure}
  \begin{center}
    \includegraphics[height=5cm, width=15cm]{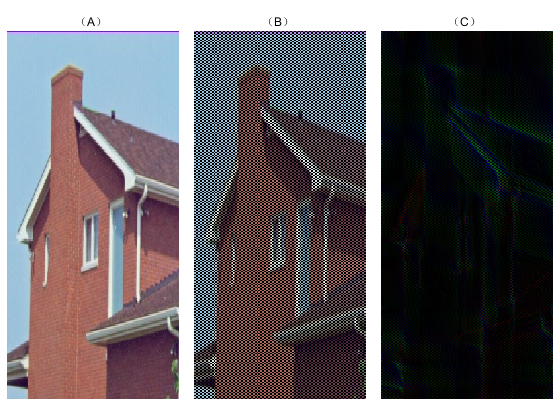}
     \caption{ The left figure (A) is the original picture, the middle figure (B) is the outcome of the lowpass filter $H_1$,the right figure (C) is the outcome of the lowpass filter $H_2$.  }
 \label{Fig333}
   \end{center}
\end{figure}
From the Figure \ref{Fig333} and the values of the SNR and PSNR, we see that      the smaller cut-off frequency,   the  lower quality of the image.
\end{example}

\begin{example}
To reduce noise, we apply the lowpass filter $H_1$ in Eq.(\ref{H1}). Gaussian noise is added to the image of the house. The figures (a) and (c) in Figure \ref{Fig444} are noised images caused by Gaussian noises, while the figures (b) and (d) are lowpass filtered images. 
The noise performance is measured in terms of the PSNR values. As shown in Table \ref{tab3}, the PSNR values of figures (b) and (d) are higher than figures (a) and (c), indicating that the lowpass filter $H_1$ can reduce noise.
  \begin{table}[h!]
\caption{PSNR comparison values for the house image which is noised by Gaussnian noises }\label{tab3}
\centering
\begin{tabular}{|c|c|}
	\hline
	$A_1$=(0,2,-1/2,2);
	$A_2$=(0,2,-1/2,4) &  PSNR \\
	\hline
	Image $(a)$ SNR=3.3530 &  6.6950\\
	 Image $(b)$ & 6.8413 \\
	\hline
	Image $(c) $  SNR = 3.4335 &  6.8104 \\
   Image $(d)$ & 6.8691\\
	\hline	
\end{tabular}
\end{table}

\begin{figure}
	\begin{center}
		\includegraphics[height=5cm,width=10cm]{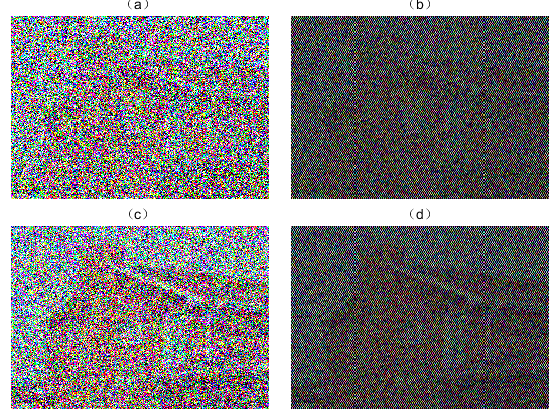}
		\caption{ The  figures $(a)$ and $(c)$ are the noised images, the  figures  $(b)$ and $(d)$ are  the outcome of the lowpass filter $H_1$ }
		\label{Fig444}
	\end{center}
\end{figure}	
\end{example}

\section{Conclusion}\label{covsec5}
In this paper, we study two novel types of convolution operators for the QLCT, namely spatial and spectral convolution operators, respectively. They are distinct in the quaternion space and are consistent once in complex or real space. The associated convolution theorems for the QLCT are derived. The spectral and spatial convolutions between two quaternionic functions can be implemented by the product of their QLCTs and the sum of their OLCTs, respectively. Furthermore, we demonstrate that the convolution theorems in the 2-D QFT
and QFrFT domains can be regarded as two
special cases of our achieved results. Four aspects of applications in convolutions are analyzed. Firstly, the correlation operations of the QLCT are developed. Secondly, the  Fredholm integral equation of the first kind involving special kernels can be solved. Thirdly, some systems of second order partial differential equations which can be transformed into the second order quaternion partial differential equations can also be solved. Finally, it is convenient to design the multiplicative filters in the QLCT domain. 
In view of the QLCT's attractive properties, its discrete version needs to be developed. In addition,  its convolution theorems should be explored in the context of quaternion random signal processing. These problems will be considered in our follow-up studies.

\section*{Acknowledgements}

Xiaoxiao Hu was supported by the Research Development Foundation of Wenzhou Medical University (QTJ18012), Wenzhou Science and Technology Bureau (G2020031), Scientific Research Task of Department of Education of Zhejiang Province (Y202147071). Dong Cheng was supported by  Guangdong Basic and Applied Basic Research Foundation (No.2019A1515111185). Kit Ian Kou was supported by University of Macau (MYRG2019-00039-FST), Science and Technology Development Fund, Macao S.A.R (FDCT/0036/2021/AGJ).


\section*{Declaration of Competing Interest}
The authors declare that they have no known competing financial interests or personal relationships that could have appeared to influence the work reported in this paper.


\begin{thebibliography}{10}
 \bibitem{hamilton}
 W.~R.~Hamilton, Elements of Quaternions, Longmans Green, London, 1866.
\bibitem{wei2020convolution}
	D.~Y.~Wei, L.~Y.~Min, Convolution and multichannel sampling for the offset
	linear canonical transform and their applications, IEEE Transactions on
	Signal Processing. 67 (2019) 6009--6024.
	
	\bibitem{PWJ2000}
	W.~Pen, {F}ourier analysisi and its application, Peking University Press, 2000.
	
	\bibitem{deng2010comments}
	B.~Deng, R.~Tao, Y.~Wang, Comments on "convolution and product theorem for the
	linear canonical transform", Signal Processing Letters, IEEE. 17~(6) (2010)
	615--616.
	
	\bibitem{shi2014generalized}
	J.~Shi, X.~P.~Liu, N.~T.~Zhang, Generalized convolution and product theorems
	associated with linear canonical transform, Signal, Image and Video
	Processing. 8~(5) (2014) 967--974.
	
	\bibitem{deng2006convolution}
	B.~Deng, R.~Tao, Y.~Wang, Convolution theorems for the linear canonical
	transform and their applications, Science in China Series {F}: Information
	Sciences. 49~(5) (2006) 592--603.
	
	\bibitem{wei2012convolution}
	D.~Y.~Wei, Q.~W.~Ran, Y.~M.~Li, A convolution and correlation theorem for the
	linear canonical transform and its application, Circuits, Systems, and Signal
	Processing. 31~(1) (2012) 301--312.
	
	\bibitem{wei2009convolution}
	D.~Y.~Wei, Q.~W.~Ran, Y.~M.~Li, J.~Ma, L.~Y.~Tan, A convolution and product
	theorem for the linear canonical transform, Signal Processing Letters, IEEE.
	16~(10) (2009) 853--856.
	
	\bibitem{kou2013uncertainty}
	K.~I.~Kou, J.~Y.~Ou, J.~Morais, On uncertainty principle for quaternionic
	linear canonical transform, in: Abstract and Applied Analysis., Vol. 2013,
	Hindawi Publishing Corporation, 2013.
	
	\bibitem{kou2014asymptotic}
	K.~I.~Kou, J.~Morais, Asymptotic behaviour of the quaternion linear canonical
	transform and the Bochner--Minlos theorem, Applied Mathematics and
	Computation. 247 (2014) 675--688.
	
	\bibitem{wang2021robust}
	Y.~L.~Wang, K.~I.~Kou, C.~M.~Zou, Y.~Y.~Tang, Robust sparse representation in
	quaternion space, IEEE Transactions on Image Processing. 30 (2021) 3637--3649.
	
	\bibitem{ell2013quaternion}
	T.~A.~Ell, Quaternion fourier transform: re-tooling image and signal processing
	analysis, in: Quaternion and Clifford {F}ourier  Transforms and Wavelets.,
	Springer, 2013, pp. 3--14.
	
	\bibitem{le2014instantaneous}
	N.~Le~Bihan, S.~J.~Sangwine, T.~A.~Ell, Instantaneous frequency and amplitude
	of orthocomplex modulated signals based on quaternion Fourier transform,
	Signal Processing. 94 (2014) 308--318.
	
	\bibitem{bahri2013convolution}
	M.~Bahri, R.~Ashino, R.~Vaillancourt, Convolution theorems for quaternion
	fourier transform: properties and applications, in: Abstract and Applied
	Analysis., Vol. 2013, Hindawi Publishing Corporation, 2013.
	\bibitem{bujack2014convolution}
     R.~Bujack, H.~De.~Bie, N.~De.~Schepper, G.~Scheuermann, Convolution products for hypercomplex Fourier transforms, Journal of Mathematical Imaging and Vision. (2014) 48:606–624.
	
	\bibitem{hitzer2016general}
	E.~Hitzer, General two-sided quaternion {Fourier} transform, convolution and
	mustard convolution, Advances in Applied Clifford Algebras. (2016) 1--15.
	
	\bibitem{de2015connecting}
	H.~De~Bie, N.~De~Schepper, T.~A.~Ell, K.~Rubrecht, S.~J.~Sangwine, Connecting
	spatial and frequency domains for the quaternion {Fourier} transform, Applied
	Mathematics and Computation. 271 (2015) 581--593.
	
	\bibitem{guanlei2008fractional}
	G.~L.~Xu, X.~T.~Wang, X.~G.~Xu, Fractional quaternion {Fourier} transform,
	convolution and correlation, Signal Processing. 88~(10) (2008) 2511--2517.
	
	\bibitem{serbes2010optimum} 
	A.~Serbes, L.~Durak, Two-sided fractional quaternion {Fourier} transform and its
	application, Journal of Inequalities and Applications. 1 (2021) 121.
	
	\bibitem{bhat2022wvd}
	M.~Y.~Bhat, A.~H.~Dar, Convolution and correlation theorems for {Wigner-Ville}
	distribution associated with the quaternion offset linear canonical
	transform, Signal, Image and Video Processing. 16 (2022) 1235--1242.
	
	\bibitem{bhat2022algebra}
	M.~Y.~Bhat, A.~H.~Dar, The algebra of $2D$ Gabor quaternion offset linear
	canonical transform and uncertainty principles, The Journal of Analysis. 30
	(2022) 637--649.
	
	\bibitem{guo2011reduced}
	L.~Q.~Guo, M.~Zhu, X.~H.~Ge, Reduced biquaternion canonical transform,
	convolution and correlation, Signal Processing. 91~(8) (2011) 2147--2153.
	
	\bibitem{hu2021sampling}
	X.~X.~Hu, D.~Cheng, K.~I.~Kou, Sampling formulas for $2D$  quaternionic signals
	associated with various quaternion {Fourier} and linear canonical transforms,
	Frontiers of Information Technology and Electronic Engineering. 23~(3) (2022)
	463--478.
	
	\bibitem{chen2015pitt}
	L.~P.~Chen, K.~I.~Kou, M.~S.~Liu, Pitt's inequality and the uncertainty
	principle associated with the quaternion Fourier transform, Journal of
	Mathematical Analysis and Applications. 423~(1) (2015) 681--700.
	
	\bibitem{yang2014uncertainty}
	Y.~Yang, K.~I.~Kou, Uncertainty principles for hypercomplex signals in the
	linear canonical transform domains, Signal Processing. 95 (2014) 67--75.
	
	\bibitem{kou2016envelope}
	K.~I.~Kou, M.~S.~Liu, J.~P.~Morais, C.~M.~Zou, Envelope detection using
	generalized analytic signal in {2D} {QLCT} domains, Multidimensional Systems and
	Signal Processing. (2016) 1--24.
	
	\bibitem{hu2016quaternion}
	X.~X.~Hu, K.~I.~Kou, Quaternion {Fourier} and linear canonical inversion
	theorems, Mathematical Methods in the Applied Sciences. 40 (2017) 2421--2440.
	
	\bibitem{kou2019Plancherel}
	K.~I.~Kou, M.~S.~Liu, C.~M.~Zou, Plancherel theorems of quaternion {Hilbert}
	transforms associated with linear canonical transforms, Advances in Applied
	Clifford Algebras. 30~(2019)9.
	
	\bibitem{hu2022sampling}
	X.~X.~Hu, K.~I.~Kou, Sampling formulas for non-bandlimited quaternionic
	signals, Signal, Image and Video Processing. (2022).
	{\path{doi:10.1007/s11760-021-02110-1}}.
	
	\bibitem{sam2021new}
	S.~Saima, B.~Z.~Li, Quaternionic one-dimensional linear canonical transform,
	Optik. 244 (2021) 166914.
	
	\bibitem{bahri2019Two}
	M.~Bahri, R.~Ashino, Two-dimensional quaternion linear canonical transform:
	Properties, convolution, correlation, and uncertainty principle, Journal of
	Mathematics. (2019)1062979.
	
	\bibitem{li2021new}
	Z.~W.~Li, W.~B.~Gao, B.~Z.~Li, A new kind of convolution, correlation and
	product theorems related to quaternion linear canonical transform, Signal,
	Image and Video Processing. 15~(1) (2021) 103--110.
	
	\bibitem{ell2014quaternion}
	T.~A.~Ell, N.~Le~Bihan, S.~J.~Sangwine, Quaternion Fourier transforms for
	signal and image processing, John Wiley \& Sons, 2014.
	
	\bibitem{wei2012new}
	D.~Y.~Wei, Q.~Ran, Y.~Li, New convolution theorem for the linear canonical
	transform and its translation invariance property, Optik-International
	Journal for Light and Electron Optics. 123~(16) (2012) 1478--1481.
	
	\bibitem{sangwine1999hypercomplex}
	S.~J.~Sangwine, T.~A.~Ell, Hypercomplex auto-and cross-correlation of color
	images, in: Image Processing, 1999. ICIP 99. Proceedings. 1999 International
	Conference on, Vol.~4, IEEE, 1999, pp. 319--322.
	
	\bibitem{ell2000hypercomplex}
	T.~A.~Ell, S.~J.~Sangwine, Hypercomplex Wiener-Khintchine theorem with
	application to color image correlation, in: 2000 International Conference on
	Image Processing, 2000. Proceedings., Vol.~2, IEEE, 2000, pp. 792--795.
	
	\bibitem{took2011augmented}
	C.~C.~Took, D.~P.~Mandic, Augmented second-order statistics of quaternion
	random signals, Signal Processing. 91~(2) (2011) 214--224.
	
	\bibitem{navarro2016semi}
	J.~Navarro~Moreno, J.~C. Ruiz-Molina, Semi-widely linear estimation of
	c$\eta$-proper quaternion random signal vectors under gaussian and stationary
	conditions, Signal Processing. 119 (2016) 56--66.
	
	\bibitem{ginzberg2013quaternion}
	P.~Ginzberg, A.~T.~Walden, Quaternion var modelling and estimation, IEEE
	transactions on signal processing. 61~(1-4) (2013) 154--158.
	
	\bibitem{gou2015three}
	X.~M.~Gou, Z.~W.~Liu, W.~Liu, Y.~G.~Xu, Three-dimensional wind profile
	prediction with trinion-valued adaptive algorithms, in: 2015 IEEE
	International Conference on Digital Signal Processing (DSP)., IEEE, 2015, pp.
	566--569.
	
	\bibitem{kress1989linear}
	R.~Kress, V.~Maz'ya, V.~Kozlov, Linear integral equations, Vol.~82, Springer,
	1989.
	
	\bibitem{chan1996conjugate}
	R.~H.~Chan, M.~K.~Ng, Conjugate gradient methods for toeplitz systems, SIAM
	review. 38~(3) (1996) 427--482.
	
	\bibitem{rajan2003simultaneous}
	D.~Rajan, S.~Chaudhuri, Simultaneous estimation of super-resolved scene and
	depth map from low resolution defocused observations, IEEE Transactions on
	Pattern Analysis and Machine Intelligence. 25~(9) (2003) 1102--1117.
	
	\bibitem{polyanin2008handbook}
	A.~D.~Polyanin, A.~V.~Manzhirov, Handbook of integral equations, CRC press,
	2008.
	
	\bibitem{florian2001functional}
	H.~Florian, N.~Ortner, F.~Schnitzer, W.~Tutschke, Functional-analytic and
	complex methods, their interactions, and applications to partial differential
	equations, in: Proceedings of the Workshop held at Graz University of
	Technology, World Scientific, 2001.
	
	\bibitem{smirnov1964course}
	V.~I.~Smirnov, D.~BROWN, A Course of Higher Mathematics Translated (from the
	16th Russian Edition) by DE Brown. Translation Edited and Additions Made by
	IN Sneddon, Etc, Pergamon Press, 1964.
	
	\bibitem{pei2001efficient}
	S.~C.~Pei, J.~J.~Ding, J.~H.~Chang, Efficient implementation of quaternion
{Fourier} transform, convolution, and correlation by 2-D complex FFT, IEEE
	Transactions on Signal Processing. 49~(11) (2001) 2783--2797.
	
\end{thebibliography}

\end{document}